\documentclass[12pt]{amsart}
\usepackage{amsmath}
\usepackage{amsthm}
\usepackage{amsfonts}
\usepackage{mathdots}

\usepackage[T1]{fontenc}
\usepackage[utf8]{inputenc}

\usepackage{times}

\usepackage{amssymb}
\usepackage{hyperref}

\author[S.~Belinschi]{Serban T.~Belinschi}
\address{CNRS,  
Institut de Math\'ematiques de Toulouse, and Department of Mathematics and Statistics,
Queen's University.
Institut de Math\'ematiques de Toulouse,
Universit\'e Paul Sabatier,
118, route de Narbonne
F-31062 Toulouse Cedex 9, FRANCE} \email{serban.belinschi@math.univ-toulouse.fr}

\author[P.~\'Sniady]{Piotr \'Sniady}
\address{
Institute of Mathematics, Polish Academy of Sciences,
\mbox{ul.~\'Sniadec\-kich 8,} \linebreak 00-956 Warszawa, Poland
} 
\email{psniady@impan.pl}

\author[R.~Speicher]{Roland Speicher}
\address{Universit\"at des Saarlandes, Fachrichtung 6.1 - Mathematik, Postfach
151150, 66041 Saarbr\"ucken, Germany} 
\email{speicher@math.uni-sb.de}

\keywords{eigenvalues of non-Hermitian random matrices, Brown measure, non-normal operators}
\subjclass{15B52 (Primary), 46L54, 46L10 (Secondary)}

\title[Brown measure of polynomials in free variables]
{Eigenvalues of non--Hermitian random matrices and Brown measure of non-normal operators: 
Hermitian reduction and linearization method}

\theoremstyle{plain}
\newtheorem{lemma}{Lemma}
\newtheorem{theorem}[lemma]{Theorem}
\newtheorem{proposition}[lemma]{Proposition}

\newtheorem*{conjecture}{Conjecture}

\theoremstyle{definition}

\theoremstyle{remark}
\newtheorem*{remark}{Remark}

\newcommand{\x}{X}
\newcommand{\G}{{\mathbf{G}}}
\newcommand{\F}{{\mathbf{F}}}

\newcommand{\A}{{\mathfrak{A}}}
\newcommand{\B}{{\mathfrak{B}}}

\newcommand{\E}{{\mathbb{E}}}
\newcommand{\C}{{\mathbb{C}}}
\newcommand{\R}{{\mathbb{R}}}

\newcommand{\El}{{\mathcal{L}}}
\newcommand{\M}{{\mathcal{M}}}

\newcommand{\gwia}{^{\star}}

\newcommand{\ddd}{d}

\DeclareMathOperator{\supp}{supp}

\DeclareMathOperator{\tr}{tr} \DeclareMathOperator{\Tr}{Tr}

\begin{document}

\begin{abstract}
We study the Brown measure of certain non--Hermitian
operators arising from Voiculescu's free probability theory. 
Usually those operators appear as the limit in $\star$-moments
of certain ensembles of non--Hermitian random matrices, and the 
Brown measure gives then a canonical candidate for the limit eigenvalue
distribution of the random matrices.
A prominent class for our operators is given by polynomials in $\star$-free variables.
Other explicit 
examples include $R$--diagonal elements and elliptic elements, for
which the Brown measure was already known, and a new class of
triangular--elliptic elements. 
Our method for the calculation of the Brown measure is
based on a rigorous mathematical treatment of the Hermitian
reduction method, as considered in the physical literature, combined with subordination and the linearization trick. 
\end{abstract}

\maketitle

\section{Introduction}

\subsection{Eigenvalues of non-Hermitian random matrices}

The study of the eigenvalue distribution  of non-Hermitian random
matrices is regarded as an important and interesting problem,
especially in the mathematical physics literature. Unfortunately,
most of the methods used for the study of Hermitian random
matrices fail in the non-Hermitian case, which makes the latter
very difficult.

\subsection{Convergence of $\star$--moments. Free probability theory}
\label{subsec:convergence} Usually we are interested in the
behavior of the random matrix eigenvalues in the limit as the size
of the matrix tends to infinity. It is therefore natural to ask:
does a given sequence of random matrices converge in one sense or
another to some (infinite-dimensional) object as the size of
the matrix tends to infinity? It would be very tempting to study
this limit instead of the sequence of random matrices itself.

In order to perform this program, we will use the notion of a
$W\gwia$-probability space (which is a von Neumann algebra $\A$
equipped with a tracial, faithful, normal state
$\phi:\A\rightarrow\C$). The algebra $\A_N=\El^{\infty-}(\Omega,\M_N)$
of $N\times N$ random matrices with all moments finite equipped
with a tracial state $\phi_N(x)= \frac{1}{N}\E \Tr x$ for
$x\in\A_N$ fits well into this framework (to be very precise: the
definition of a von Neumann algebra requires its elements to be
bounded which is not the case for the most interesting examples of
random matrices, but this small abuse of notation will not cause
any problems in the following).

We say that a sequence  of random matrices $(A_N)$, where
$A_N\in\A_N$, converges in $\star$--moments to some element
$x\in\A$ if for every choice of $s_1,\dots,s_n\in\{1,\star\}$ we
have
$$ \lim_{N\to\infty} \phi_N(A_N^{s_1} \cdots A_N^{s_n})=
\phi(x^{s_1} \cdots x^{s_n}).$$ It turns out that many classes of
random matrices have a limit in a sense of $\star$--moments and
the limit operator can be found by the means of free probability
theory \cite{VDN,HP,MSbook}.

\subsection{Brown measure}
The Brown measure \cite{Brown} is an analogue of the density
of eigenvalues for elements of $W\gwia$--probability spaces. Its
great advantage is that it is well--defined not only for
self-adjoint or normal operators; furthermore for random matrices
it coincides with the mean empirical eigenvalue distribution. We
recall the exact definition in Section \ref{subsec:brown}.

\subsection{The main tools: $\mathbf{(i)}$ Hermitian reduction method}
In this article we study rigorously the idea of Janik, Nowak, Papp
and Zahed \cite{JanikNowakPappZahed1997} which was later used
in the papers
\cite{FeinbergZeeNongaussian,FeinbergZeeHermitianreduction,
FeinbergScalettarZee} under the name \emph{Hermitian reduction
method}.

As we shall see in Section \ref{subsec:brown}, the Brown
measure $\mu_x$ of an element $x$ is closely related to the Cauchy
transform $G_x(\lambda)=\phi\big( (x-\lambda)^{-1} \big)$ of $x$.
The asymptotic expansion of $G_x(\lambda)$ at infinity is given by
the sequence of moments of $x$ and for this reason it can be
computed explicitly by the means of free probability theory
for many operators $x$.
However, given a sequence of moments $\{m_n=\int t^n\,d\mu_x(t)\}_{n\in\mathbb N}$, there 
usually are many probability measures supported in $\mathbb C$ which have $\{m_n\}_n$ as
their sequence of moments. In particular, 
if $x$ is not Hermitian
then 
the series $G_x(\lambda)=\sum_{n=0}^\infty\phi(x^n)\lambda^{-n-1}$ alone does not determine the distribution of $x$.

The idea is to arrange the operator $x$ and its Hermitian conjugate $x\gwia$ into a
$2\times 2$ matrix and to consider a matrix--valued Cauchy
transform
$$ \mathbf{G}_{\epsilon}(\lambda)=
\phi \left( \begin{bmatrix} i \epsilon & \lambda-x \\
\overline{\lambda}-x\gwia & i \epsilon
\end{bmatrix}^{-1} \right)\in\M_2(\C)$$
which depends on an additional parameter $\epsilon>0$. The above function makes sense for all such $\epsilon>0$; it will be viewed as a restriction of an analytic map defined 
on a matricial upper half-plane. 
As we shall see, the limit $\epsilon\to 0$
gives us access to the original Cauchy transform $G_x(\lambda)$
and therefore to the Brown measure of $x$.

This method appears to be extremely simple and indeed its
applications in the physics literature involved very
short and simple calculations. However, from a mathematical point
of view they are often far from being rigorous. In this article we
would like to put the Hermitian reduction method on solid ground.

The main ingredient in our approach is a recent progress in \cite{BMS} on
the analytic description of operator-valued free convolutions, relying on
the idea of subordination. Roughly speaking, subordination usually yields
an analytic description of the relevant equations (say, for the operator-valued
Cauchy transforms) which are not only valid in some neighborhood of infinity, but everywhere in the (operator-valued) complex upper half plane. Since
recovering the desired distribution relies on the knowledge of the
Cauchy transform close to the real axis the analytic description everywhere in the complex upper half plane is crucial for a rigorous treatment. As examples for such explicit calculations we will present a detailed analysis of two interesting classes of
non--Hermitian random matrices and the corresponding
non--Hermitian operators (namely, so-called $R$-diagonal and elliptic-triangular operators). 

\subsection{The main tools: $\mathbf{(ii)}$ linearization method}
It seems that by mimicking the methods
presented in this article it should be possible to calculate the
Brown measure of virtually any operator described in terms of
free probability. As a very general class of such operators we will in
particular consider the problem of arbitrary (in general, non-self-adjoint)
polynomials in free variables. 
In the above mentioned paper \cite{BMS} the corresponding problem for self-adjoint polynomials in self-adjoint free variables was solved by invoking what is known as \emph{linearization trick}, which allows us
to reduce the polynomial problem to an operator-valued linear problem. 
The same reduction works in the non-self-adjoint case, and we will show how the combination of the Hermitization and the linearization methods will result in 
an algorithm for calculating the Brown measure of a polynomial in $\star$-free variables out of the $\star$-distributions of its variables; in particular, if the variables are normal then this gives a way to calculate the Brown  measure of the polynomial out of the Brown measures of the variables.

We will then also start an investigation on the qualitative features of the simplest polynomial, namely the sum of two $\star$-free operators.  We will, in particular, address the question how eigenvalues in the sum can arise from eigenvalues in summands. This is related (but not equivalent) to the question
of atoms in the corresponding Brown measures. 

\subsection{Overview of this article and statement of results}
In Section \ref{sec:preliminaries} we recall the definitions of
the Brown measure and the Cauchy transform as well as their
regularized versions and their connection with the Hermitian
reduction method. We also recall some basic tools of Voiculescu's
free probability theory.

In Section \ref{section:linearization} we present the linearization trick
and show how it combines with the Hermitian reduction method to yield an
algorithm
for the computation of the Brown measure of polynomials in free variables.

In Section \ref{section:polynomials} we will then address in more detail
polynomials in free variables; in particular, we present analytic properties
of the Brown measure of the sum of two $\star$-free variables.

The next two sections will then deal with important special classes of
non-normal operators and we will (re)derive explicit expressions for their Brown measures. 

In Section \ref{sec:rdiagonal} we study the class of $R$--diagonal
operators which are limits of bi-unitarily invariant
random matrices. We calculate by our methods the Brown measure of $R$--diagonal
operators and rederive in this way the results of Haagerup and Larsen \cite{HaagerupLarsen}. That this Brown measure is actually the limit of the eigenvalue distributions of the corresponding random matrix models was recently proved by Guionnet, Krishnapur, and Zeitouni \cite{GZ}.
We also note that results of the present paper are used in \cite{BNST} to determine
the asymptotics of the eigenvector overlap for those random matrices.

In Section \ref{sec:triangular} we study certain non--Hermitian
Gaussian random matrices the entries of which above the diagonal,
informally speaking, have a different covariance than the entries
below the diagonal. We describe the operators which arise as
limits of such random matrices, and we use the Hermitian reduction
method to calculate their Brown measures. 
By a result of \'Sniady \cite{Sn} we know that this agrees in this case with the limit
distribution of the eigenvalues of these random matrices.

In Section \ref{sec:discontinuity} we then briefly discuss the problem of
discontinuity of the Brown measure. Roughly speaking, the
eigenvalues of non--Hermitian matrices do not depend on the matrix
$\star$-moments in a continuous way and therefore the Brown measure of
some operator might not be related to the eigenvalues of matrices
which converge in $\star$-moments to this operator. However, one expects that
for natural choices of random matrices such a convergence should hold.

\section{Preliminaries}
\label{sec:preliminaries}

\subsection{Cauchy transform and Brown measure}
\label{subsec:brown}

\subsubsection{Cauchy transform} Let $\mu$ be a probability measure
on the complex plane $\C$. We define its Cauchy transform 
as the analytic function 
\begin{equation}
\label{eq:Cauchy} G_\mu(\lambda)= \int_{\C} \frac{1}{\lambda-z}
d\mu(z)
\end{equation}
for $\lambda \notin \supp \mu$.
It is known \cite{Garnett} 
that the integral $\int_{\C} \frac{1}{\lambda-z}\, d\mu(z)$ is in fact
well-defined everywhere outside a set of $\mathbb R^2$-Lebesgue measure zero, and thus
$G_\mu$ will be viewed from now on as a function on all of $\mathbb C$, whose analyticity
will, however, be guaranteed only outside the support of $\mu$.
The measure $\mu$ can be extracted from its
Cauchy transform by the formula
\begin{equation}
\label{eq:recover} \mu = \frac{1}{\pi} \frac{\partial}{\partial
\bar{\lambda}} G_{\mu}(\lambda),
\end{equation}
where as usual
$$\frac{\partial}{\partial \lambda}= \frac{1}{2} \left(
\frac{\partial}{\partial (\Re \lambda)} - i
\frac{\partial}{\partial (\Im \lambda)} \right), \qquad
\frac{\partial}{\partial \bar{\lambda}}= \frac{1}{2} \left(
\frac{\partial}{\partial (\Re \lambda)} + i
\frac{\partial}{\partial (\Im \lambda)} \right)$$ denote the
derivatives in the Schwartz distribution sense; 
thus, \eqref{eq:recover} should be understood in the distributional sense too.

\subsubsection{Spectral measure of self-adjoint operators}
\label{subseubsec:example1} Let $\A$ be a von Neumann algebra
equipped with a normal faithful tracial state $\phi$. Every
self-adjoint operator $x\in\A$ can be written as a spectral
integral
$$x= \int_{\R} \lambda\ dE(\lambda),$$
where $E$ denotes the operator--valued spectral measure of $x$. It
is natural to consider a probability measure $\mu_x$ on $\C$ given
by
\begin{equation}
\label{eq:Brown01} \mu_x(Z) = \phi\big( E(Z) \big)
\end{equation}
for any Borel set $Z\subseteq\C$.

In a full analogy with \eqref{eq:Cauchy} we consider the Cauchy
transform of $x$ given by
\begin{equation}
\label{eq:Cauchy2} G_{x}(\lambda)=\phi\left( (\lambda-x)^{-1}
\right).
\end{equation}
Then the spectral measure $\mu_x$ as defined by \eqref{eq:Brown01} %{eq:Cauchy2}
can be recovered from \eqref{eq:Cauchy2} via \eqref{eq:recover}. 

\subsubsection{Spectral measure of matrices}
\label{subseubsec:example2} Let $\A=\M_N$ be the matrix algebra
equipped with the tracial state $\phi=\tr$, where $\tr x=
\frac{1}{N} \Tr x$ is a normalized trace. Let
$\lambda_1,\dots,\lambda_N\in\C$ be the eigenvalues (counted with
multiplicities) of a given matrix $x\in\M_N$. We define $\mu_x$ to
be the (normalized) counting measure of the set of eigenvalues
\begin{equation}
\label{eq:Brown02}
\mu_{x}=\frac{\delta_{\lambda_1}+\cdots+\delta_{\lambda_N}}{N}.
\end{equation}
The Cauchy transform of $x$ defined by \eqref{eq:Cauchy2} is
well--defined on the set
$\C\setminus\{\lambda_1,\dots,\lambda_N\}$ and it again satisfies
\eqref{eq:recover}.

\subsubsection{Empirical eigenvalue distribution}
\label{subseubsec:example3} Let $\A=\El^{\infty-}(\Omega,\M_N)$ be
the algebra of random matrices having all moments finite. We equip
it with a state $\phi(x)=\E \tr x$. For $x\in\A$ consider a random
variable $\Omega\ni\omega\mapsto\mu_{x(\omega)}$, called empirical
eigenvalue distribution, the values of which are probability
measures on $\C$ (where $\mu_{x(\omega)}$ is to be understood as
in Section \ref{subseubsec:example2}). We define the mean eigenvalue
distribution $\mu_x$ by \begin{equation} \label{eq:Brown3} \mu_x=
\E \mu_{\omega(x)}. \end{equation}

\subsubsection{Brown measure} Let $\A$ be a von Neumann
algebra equipped with a normal faithful tracial state $\phi$.
Inspired by the above examples we might try to define the Cauchy
transform of $x\in\A$ by the formula \eqref{eq:Cauchy2} and then
define its spectral measure $\mu_x$ by \eqref{eq:recover}.
However, in general this is not possible because formula \eqref{eq:Cauchy2} requires
$\lambda$ to lie outside of the spectrum of $x$. In the non-Hermitian case, the spectrum might be a large set; in fact, it can be an arbitrary compact subset of $\mathbb C$. For such arbitrary subsets of $\mathbb C$, the moment
problem is not well-defined.
Thus, knowing the Cauchy
transform on the resolvent set of $x$ might very well not be sufficient. 
For this reason we need a more elaborate
definition of the spectral measure.

The Fuglede--Kadison determinant $\Delta(x)$ of $x\in\A$ is defined in
\cite{FugledeKadison} by
$$\log \Delta(x) = \frac{1}{2} \phi\big( \log (x x\gwia ) \big). $$
If $x$ is not invertible, the above definition should be
understood as $\Delta(x)=\lim_{\epsilon\to 0}
\Delta_{\epsilon}(x)$, where $\Delta_\epsilon$ denotes the
regularized Fuglede--Kadison determinant
$$\log \Delta_{\epsilon}(x) = \frac{1}{2} \phi\big( \log (x x\gwia + \epsilon^2) \big) $$
for $\epsilon>0$.

The Brown measure of $x\in\A$ is defined in \cite{Brown} by
\begin{equation}
\label{eq:brown} \mu_x = \frac{1}{2\pi}
\left(\frac{\partial^2}{\partial (\Re \lambda)^2}+
\frac{\partial^2}{\partial (\Im \lambda)^2} \right) \log
\Delta(x-\lambda) = \frac{2}{\pi} \frac{\partial}{\partial
\lambda}  \frac{\partial}{\partial \bar{\lambda}}  \log
\Delta(x-\lambda).
\end{equation}
The Brown measure $\mu_x$ as defined in \eqref{eq:brown} could be a
priori a Schwartz distribution but one can show that the map
$\lambda\mapsto \log \Delta(x-\lambda)$ is subharmonic and hence
$\mu_x$ is a positive measure on $\C$. In fact $\mu_x$ is a
probability measure supported on a subset of the spectrum of
$x$. One can show that for the examples from Sections
\ref{subseubsec:example1}--\ref{subseubsec:example3} the above
definition gives the correct values for \eqref{eq:Brown01},
\eqref{eq:Brown02}, and \eqref{eq:Brown3}.

Following \eqref{eq:Cauchy} we 
 define the
Cauchy transform of $x$ as
\begin{equation} G_x(\lambda)=
\int_{\C} \frac{1}{\lambda-z} d\mu_x(z). \label{eq:Cauchy3}
\end{equation}

\subsubsection{Regularized Cauchy transform and regularized Brown
measure} \label{subsubsec:regularized} For every $\epsilon>0$ the
regularized Cauchy transform
\begin{equation}
\label{eq:regular} G_{\epsilon,x}(\lambda)=  \phi\left({(\lambda-x)\gwia}
{\left((\lambda-x) (\lambda-x)\gwia+\epsilon^2\right)^{-1}}
\right)
\end{equation}
is well--defined for every $\lambda\in\C$, 
but is not an analytic function. 
In fact it was shown by Larsen
\cite{LarsenPhD} (see also \cite[Lemma 4.2]{AH}) 
and can be verified through direct arithmetic
(here it is essential to remember that $\phi$ is tracial!), that
$$ G_{\epsilon,x}(\lambda)= 2 \frac{\partial}{\partial \lambda} \log \Delta_{\epsilon}(x-\lambda).
$$
The function $\lambda\mapsto\log \Delta_{\epsilon}(x-\lambda)$
is subharmonic,
hence the regularized Brown measure defined by
\begin{equation}
\label{eq:guga} \mu_{\epsilon,x}= \frac{1}{\pi}
\frac{\partial}{\partial \bar{\lambda}} G_{\epsilon,x}(\lambda)=
\frac{2}{\pi} \frac{\partial}{\partial
\bar{\lambda}}\frac{\partial}{\partial \lambda} \log
\Delta_{\epsilon}(x-\lambda)
\end{equation}
is a positive measure on the complex plane. Integration by parts
shows that for $\epsilon\to 0$ the regularized Brown measure
$\mu_{\epsilon,x}$ converges (in the weak topology of probability
measures) towards the Brown measure
$\mu_x$ as defined by \eqref{eq:brown}. 
The comparison of \eqref{eq:guga} and \eqref{eq:recover} is a
heuristic argument that the definition of the Brown measure is a
reasonable extension of the cases from Sections
\ref{subseubsec:example1}--\ref{subseubsec:example3}.

We should probably point the reader to the fact that the measure $\mu_{\epsilon,x}$ has full
support. Indeed, if $\lambda=u+iv$, then, 
\begin{eqnarray*}
\lefteqn{
\left(\frac{\partial^2}{\partial u^2}+\frac{\partial^2}{\partial v^2}\right)\log\Delta_\epsilon
(x-\lambda)=}\\
& & \mbox{}\phi\left(\left[2-(\lambda-x)\left((\lambda-x) (\lambda-x)\gwia+\epsilon^2\right)^{-1}(\lambda-x)\gwia\right.\right.\\
& & \mbox{}-\left.(\lambda-x)\gwia\left((\lambda-x)
(\lambda-x)\gwia+\epsilon^2\right)^{-1}(\lambda-x)\right]\\
& & \mbox{}\times\left.\left((\lambda-x) (\lambda-x)\gwia+\epsilon^2\right)^{-1}\right).
\end{eqnarray*}
For $\epsilon=0$ and $\lambda$ outside the spectrum of $x$, it follows straightforwardly
that the term on the third row above (second row in the expression on the right of $=$) is
equal to one. Grouping the first row with the last and applying traciality of $\phi$ allows
us to conclude that $\left(\frac{\partial^2}{\partial u^2}+\frac{\partial^2}{\partial v^2}\right)\log\Delta(x-\lambda)=0$ for $\lambda-x$ invertible. But it makes equally
clear that, since in general $x\gwia(xx\gwia+\epsilon^2)^{-1}x\le1,$ with no equality
for $\epsilon\neq0$, the above is always positive when $\epsilon>0$ (we have used
here also the faithfulness of the trace).

Equation \eqref{eq:guga} 
implies that 
\begin{equation} G_{\epsilon,x}(\lambda)=
\int_{\C} \frac{1}{\lambda-z} d\mu_{\epsilon,x}(z),
\label{eq:Cauchy4}
\end{equation}
with the equality defined a priori only $\mathbb R^2$-Lebesgue almost everywhere, but extended by continuity
to all $\lambda\in\mathbb C$.
Since for $\epsilon\to 0$ the measures $\mu_{\epsilon,x}$ converge 
weakly to $\mu_x$, \eqref{eq:Cauchy3} and \eqref{eq:Cauchy4}
imply that the regularized Cauchy transforms $G_{\epsilon,x}$
converge to $G_{x}$ in the local $\El^1$ norms; in particular
$G_{\epsilon,x}(\lambda)\to G_{x}(\lambda)$ for almost all
$\lambda\in\C$. It should be mentioned that in fact the limit 
\begin{equation}\label{finiteness}
\lim_{\epsilon\to0}G_{\epsilon,x}(\lambda)=G_{x}(\lambda)\in{\mathbb C}
\end{equation}
exists for all $\lambda\in\mathbb C$ for which 
$$\lim_{\epsilon\to0}\phi\left(
\left((\lambda-x) (\lambda-x)\gwia+\epsilon^2\right)^{-1}\right)<\infty$$ 
(this limit
always exists and is strictly positive, unless $x$ is a multiple of the identity, but may very well be 
infinite).  
Unfortunately,  finiteness of the limit 
can usually only be guaranteed for $\lambda$ outside the spectrum of $x$.
Indeed, since 
$$G_{\epsilon,x}(\lambda)=\phi\left((\lambda-x)\gwia
\left((\lambda-x) (\lambda-x)\gwia+\epsilon^2\right)^{-1}\right),$$ 
we consider the
decomposition of $(\lambda-x)\gwia$ into four positive operators. For any operator $v\ge0$, 
\begin{eqnarray*}
0&\le&\left((\lambda-x) (\lambda-x)\gwia+\epsilon^2\right)^{-1/2}v
\left((\lambda-x) (\lambda-x)\gwia+\epsilon^2\right)^{-1/2}\\
&\le&\|v\|\left((\lambda-x) (\lambda-x)\gwia+\epsilon^2\right)^{-1}.
\end{eqnarray*}
Applying $\phi$ to the above inequalities and the monotonicity of the correspondence 
$\epsilon\mapsto\left((\lambda-x) (\lambda-x)\gwia+\epsilon^2\right)^{-1}$ allows us to 
conclude.

\subsection{Hermitian reduction method}\label{subsec:2.2}
Following the idea of the Janik, Nowak, Papp, Zahed
\cite{JanikNowakPappZahed1997},
for fixed $x\in\A$ 
let
\begin{equation}
\label{eq:X} \x=
\begin{bmatrix} 0 & x \\ x\gwia & 0 \end{bmatrix} \in \M_2(\A).
\end{equation}
This is trivially a self-adjoint element in $\M_2(\A).$ We equip
the algebra $\M_2(\A)$ with a positive conditional expectation
$\E:\M_2(\A)\rightarrow \M_2(\C)$ given by \begin{equation}
\label{eq:conditional} \E
\begin{bmatrix} a_{11} &
a_{12} \\ a_{21} & a_{22} \end{bmatrix} = \begin{bmatrix} \phi(a_{11}) & \phi(a_{12}) \\
\phi(a_{21}) & \phi(a_{22}) \end{bmatrix}. \end{equation}
Following
\cite{coalgebra}, we can define a fully matricial $\M_2(\C)$-valued Cauchy-Stieltjes transform:
for any $b\in\M_2(\C)$ which satisfies the condition that $\Im b:=(b-b^*)/2i>0$, the 
map
\begin{equation}\label{def-matr-valued-Cauchy}
\G_X(b)=\mathbb E\left[(b-X)^{-1}\right]
\end{equation}
is well defined and analytic on the set of elements $b$ for which $\Im b>0$.
In particular, for $\epsilon>0$ the element
$$\Lambda_{\epsilon}=  \begin{bmatrix} i\epsilon & \lambda \\
\bar{\lambda}  & i\epsilon \end{bmatrix} \in \M_2(\C) $$ 
belongs to the domain of $\G_X$, and
\begin{equation} \label{eq:ciemnosc} 
\G_{\epsilon}(\lambda)=
\G_X
(\Lambda_\varepsilon) = \begin{bmatrix}
g_{\epsilon,\lambda,11} & g_{\epsilon,\lambda,12}
\\ g_{\epsilon,\lambda,21} & g_{\epsilon,\lambda,22}
\end{bmatrix}=
\E \big( (\Lambda_\epsilon- \x)^{-1}\big) . \end{equation} Note
that the element $\Lambda_{0} - \x $ is self-adjoint, and for this
reason
$\Im\Lambda_\epsilon=\epsilon1$, making
the element $\Lambda_{\epsilon}- \x$
invertible whenever 
$\epsilon\neq0$, 
guarantees that \eqref{eq:ciemnosc} makes sense. 
One can easily check that
\begin{multline}
\label{eq:inwers} (\Lambda_\epsilon-\x)^{-1}= \\ \begin{bmatrix}
-i \epsilon \big((\lambda-x)(\lambda-x)\gwia+\epsilon^2\big)^{-1}
& (\lambda-x)\big((\lambda-x)\gwia (\lambda-x)+\epsilon^2\big)^{-1} \\
(\lambda-x)\gwia \big(
(\lambda-x)(\lambda-x)\gwia+\epsilon^2\big)^{-1}  & -i \epsilon
\big( (\lambda-x)\gwia(\lambda-x)+\epsilon^2 \big)^{-1}
\end{bmatrix}.
\end{multline}
Equations \eqref{eq:regular} and \eqref{eq:inwers} show that two
of the entries of $\G_{\epsilon}(\lambda)$ carry important
information, namely they coincide with the
regularized Cauchy transform, and its adjoint, respectively:
$$g_{\epsilon,21}(\lambda)=\overline{g_{\epsilon,12}(\lambda)}= G_{\epsilon,x}(\lambda).$$
It is therefore very tempting to ask what kind of information
is being carried by the other two entries
$g_{\epsilon,11}(\lambda)=g_{\epsilon,22}(\lambda)$. It was shown
by Janik et al.\ \cite{JanikNoerenbergNowakPappZahed} that if $x$
is a random matrix then in the limit $\epsilon\to 0$ these two entries
provide information about the interplay between the bases of the
left and the right eigenvectors. 

More generally, we record for future use that
\begin{multline}
\label{inverse}
\begin{bmatrix}
a&b-x\\
c-x\gwia&d
\end{bmatrix}^{-1}= \\
\begin{bmatrix}
-d\left[(b-x)(c-x\gwia)-ad\right]^{-1}&(b-x)\left[(c-x\gwia)(b-x)-ad\right]^{-1}\\
(c-x\gwia)\left[(b-x)(c-x\gwia)-ad\right]^{-1}&-a\left[(c-x\gwia)(b-x)-ad\right]^{-1}
\end{bmatrix}.
\end{multline}

The reader will probably find a great disconnect between the above and the
linearization procedure described below. Indeed, the correspondence described above
$(\lambda,\epsilon)\mapsto \G_X(\Lambda_\epsilon)$ has the significant
disadvantage of being profoundly non-analytic. Below we will work {\em only} with
analytic maps, hence our argument will never be $\Lambda_\epsilon$. However, as
we are free to evaluate the analytic map $f(z,w)=zw$ in the numbers $(z,\overline{z})$,
we shall take the liberty of evaluating certain analytic functions (like $\G_X$) in 
the matrix $\Lambda_\epsilon$, without viewing it as an argument in two variables.
However, for $b=\lambda,c=\overline{\lambda}$ fixed, we can consider the
analytic map $(z,w)\mapsto\G_X\left(\begin{bmatrix} z&\lambda\\
\overline{\lambda}&w\end{bmatrix}\right)$. A simple calculation shows that
the matrix $\Im\begin{bmatrix} z&\lambda\\
\overline{\lambda}&w\end{bmatrix}$ is strictly positive if and only if $\Im z>0$ and 
$\Im z\Im w>\frac{|\lambda-\overline{\overline{\lambda}}|}{4}=0$, 
i.e.~if and only if $z$ and $w$ are in the upper half plane $\mathbb C^+$ of $\mathbb C$.

Since $\G_X$ maps the matricial upper half-plane into the matricial lower half-plane,
$\F_X(b):=\G_X(b)^{-1}$ is well defined and $\Im \F_X(b)\ge\Im b$. In the case that $x$ is not a multiple of the identity, this
inequality is in fact strict. 
One can show \cite[Proposition 2.15]{BMS} %\cite{notewithTobias} actually this is in our paper{BMS}.
that the strict inequality $\Im \F_X(b)>\Im b$ can fail only when there exists a non-zero projection
$p\in\mathcal M_2(\mathbb C)$ so that $pX=p\mathbb E[X]$. Since nontrivial projections in 
$\mathcal M_2(\mathbb C)$ are necessarily of the form $$\begin{bmatrix} \alpha&j\sqrt{\alpha
(1-\alpha)}\\
\overline{j}\sqrt{\alpha
(1-\alpha)}&1-\alpha\end{bmatrix},$$ this would require that both of $\alpha x=\alpha\phi(x)$ and
$(1-\alpha)x\gwia=(1-\alpha)\overline{\phi(x)}$, as equality of operators, hold. This 
contradicts the assumption $x$ is not a multiple of the identity. 

\section{Linearization}\label{section:linearization}
\subsection{$\star$-distribution of $x$ out of $\G_{X\otimes1_n}$}\label{2.3}

In the previous section we have shown that $\G_X$ includes all the information about 
the Brown measure of $x$. But we will be interested in knowing the Brown measure of
any polynomial $P(x_1,\dots,x_k)$ in $\star$-free random variables $x_1,\dots,x_k$ in terms
of the Brown measures (maybe less than full $\star$-distributions)
of the individual random variables. In order to be able to encapsulate all that information
and efficiently manipulate it, we will use what is called the fully matricial extension \cite{coalgebra}
of $\G_X$:
$$
\G_{X\otimes1_n}(b)=(\mathbb E\otimes\mathrm{Id}_{\mathcal M_n})\left[\left(b-X\otimes1_n\right)^{-1}\right],
$$
where $b\in\mathcal M_n(\mathcal M_2(\mathbb C))$ is so that $b-X\otimes1_n$ is invertible
(in particular, this holds true if $\Im b>0$). It is known \cite{coalgebra} that
$\G_{X\otimes1_n},n\in\mathbb N$ encodes all $\mathcal M_2(\mathbb C)$-moments of $X$,
and hence all $\star$-moments of $x$.

\subsection{Linearization} The work in \cite{BMS} introduces an iterative method to compute
spectra of self-adjoint polynomials in free variables. This is based on a linearization trick
introduced in \cite{A}. 
(Versions of this linearization trick have a long history in different settings, see \cite{BR,HMV,HT}.)  In this subsection we shall show how the 
Hermitian reduction method combines with the linearization trick to allow for the
computation of the Brown measure of any polynomial in free random variables. 

To begin with, let us linearize an arbitrary monomial: assume we desire to compute the
Brown measure of $x_1x_2\cdots x_k.$ 
We will assume that any two neighbouring
elements are $\star$-free from each other. However, it is not necessary that {\em all}
$x_1,x_2,\dots ,x_k$ are free.
The Hermitian reduction requires us to build
\begin{equation}\label{monome}
\begin{bmatrix}
i\epsilon & \lambda-x_1x_2\cdots x_k\\
\overline{\lambda}-(x_1x_2\cdots x_k)^* & i\epsilon
\end{bmatrix}^{-1}.
\end{equation}
In order to be able to use the freeness of the elements involved, we will have to separate
them in sums of (possibly quite large) matrices. Observe first that
$$
\begin{bmatrix}
0 & 1\\
1 & 0
\end{bmatrix}\begin{bmatrix}
0 & 1\\
x_1 & 0
\end{bmatrix}\begin{bmatrix}
0 & x_2\\
x_2\gwia & 0
\end{bmatrix}\begin{bmatrix}
0 & x_1\gwia\\
1 & 0
\end{bmatrix}\begin{bmatrix}
0 & 1\\
1 & 0
\end{bmatrix}=\begin{bmatrix}
0 & x_1x_2\\
x_2\gwia x_1\gwia & 0
\end{bmatrix}.
$$
Thus, by induction, if we have a matrix $
\begin{bmatrix}
0 & x_2\cdots x_{k}\\
(x_2\cdots x_{k})^* & 0
\end{bmatrix}$, we will obtain $
\begin{bmatrix}
0 & x_1x_2\cdots x_k\\
(x_1x_2\cdots x_k)^* & 0
\end{bmatrix}$ as the product
$$
\begin{bmatrix}
0 & 1\\
1 & 0
\end{bmatrix}\begin{bmatrix}
0 & 1\\
x_1 & 0
\end{bmatrix}\begin{bmatrix}
0 & x_2\cdots x_k\\
(x_2\cdots x_k)^* & 0
\end{bmatrix}
\begin{bmatrix}
0 & x_1\gwia\\
1 & 0
\end{bmatrix}\begin{bmatrix}
0 & 1\\
1 & 0
\end{bmatrix}.
$$
For the linearization trick, it would however be convenient to write this product 
explicitly: if $X_j=\begin{bmatrix}
0 & x_j\\
x_j\gwia & 0
\end{bmatrix}$, $\tilde{X}_j=\begin{bmatrix}
0 & 1\\
x_j & 0
\end{bmatrix}$, and $\mathcal J=\begin{bmatrix}
0 & 1\\
1 & 0
\end{bmatrix}$, then 
\begin{multline*}
\lefteqn{\begin{bmatrix}
0 & x_1\cdots x_k\\
(x_1\cdots x_k)^* & 0
\end{bmatrix}}\\
 =  \mathcal J\tilde{X}_1\cdots\mathcal J\tilde{X}_{k-2}\mathcal J\tilde{X}_{k-1}X_k\tilde{X}_{k-1}\gwia\mathcal J\tilde{X}_{k-2}\gwia\mathcal J\cdots
\tilde{X}_1\gwia\mathcal J.
\end{multline*}
Now we shall linearize the right-hand monomial as a separate entity over $\mathcal M_2(\C)$.
For simplicity, we shall denote $Y_j=\mathcal J\tilde{X}_j$. The linearization of 
$Y_1Y_2\cdots Y_{k-1}X_kY_{k-1}\gwia\cdots Y_2\gwia Y_1\gwia$ is then performed by the
matrix
$$
\mathbb X:=\begin{bmatrix}
&  &  &  &  &  &  & Y_1\\
&  &  &  &  &  &Y_2&-1_2\\
& &   &  &  &\iddots&\iddots& \\
&  &  &  &X_k&-1_2& & \\
& & &\iddots&\iddots& & & \\
& &Y_2\gwia&-1_2& & & &\\
&Y_1\gwia&-1_2& & & & &
\end{bmatrix}.
$$
This is a $(4k-2)\times(4k-2)$ matrix, and the entries shown are $2\times2$ matrices, with
$1_2$ being the $2\times 2$ identity matrix. All empty spaces correspond to zero entries. Let
$$
b_{1}=b\otimes e_{1,1}=\begin{bmatrix}
b & 0 & \cdots & 0\\
0 & 0 & \cdots & 0\\
\vdots & \vdots & \cdots & \vdots\\
0 & 0 & \cdots & 0
\end{bmatrix}
$$
Then the matrix
$$
b_1-\mathbb X=\begin{bmatrix}
& b &  &  &  &  &  & Y_1\\
&  &  &  &  &  &Y_2&-1_2\\
& &   &  &  &\iddots&\iddots& \\
&  &  &  &X_k&-1_2& & \\
& & &\iddots&\iddots& & & \\
& &Y_2\gwia&-1_2& & & &\\
&Y_1\gwia&-1_2& & & & &
\end{bmatrix}
$$
is invertible, and 
$$
(b_1-\mathbb X)^{-1}=
\begin{bmatrix}
(b-Y_1Y_2\cdots Y_{k-1}X_kY_{k-1}\gwia\cdots Y_2\gwia Y_1\gwia)^{-1} & \ast & \cdots &
\ast\\
\ast &\ast & \cdots & \ast\\
\cdots &\cdots & \cdots & \cdots\\
\ast &\ast & \cdots & \ast
\end{bmatrix},
$$
where $\ast$ represents some unspecified entries.
Putting 
$$b=
\begin{bmatrix}
i\epsilon & \lambda\\
\overline{\lambda}& i\epsilon
\end{bmatrix}^{-1}$$ 
in the last equation and picking then from  $(b_1-\mathbb X)^{-1}$ the ${(1,1),(1,2),(2,1),(2,2)}$ entries will provide \eqref{monome}, as desired.
Note that if $P=\sum_j q_jx_{i_1,j}x_{i_2,j}\cdots x_{i_{k(j)},j}$, then 
$$
\begin{bmatrix}
0 & P\\
P\gwia & 0
\end{bmatrix}=\sum_j\begin{bmatrix}
0 & x_{i_1,j}x_{i_2,j}\cdots (q_jx_{i_{k(j)},j})\\
(x_{i_1,j}x_{i_2,j}\cdots (q_jx_{i_{k(j)},j}))\gwia & 0
\end{bmatrix}.
$$
In order to avoid having to deal with the scalars when writing the linearization matrix, we will
``merge'' them into $x_{i_{k(j)},j}$, so that now $X_{k(j)}=\begin{bmatrix}
0 & q_jx_{i_{k(j)},j}\\
\overline{q_j}x_{i_{k(j)},j}\gwia & 0
\end{bmatrix}$.  As shown in the proof of \cite[Proposition 3.4]{BMS}, the linearization of 
$\begin{bmatrix}
0 & P\\
P\gwia & 0
\end{bmatrix}$ will then be obtained simply by ``stacking'' the corresponding linearizations
of the (it should be noted!) self-adjoint $2\times 2$-matrix monomials: we shall denote $u_j=(0,
\dots,0,Y_{i_1,j})$ and $\mathbb Y_j$ the matrix 
$$
\mathbb Y_j=\begin{bmatrix}
&  &  &  &  &Y_{i_2,j}&-1_2\\
&   &  &  &\iddots&\iddots& \\
 &  &  &X_{i_{k(j)},j}&-1_2& & \\
& &\iddots&\iddots& & & \\
&Y_{i_2,j}\gwia&-1_2& & & &\\
&-1_2& & & & &
\end{bmatrix},
$$
so that $\mathbb X_j$ is obtained from $\mathbb Y_j$ by adding one first row $u_j$ and
one first column $u_j\gwia$:
$$
\mathbb X_j=\begin{bmatrix}
0& u_j\\
u_j\gwia& \mathbb Y_j
\end{bmatrix}.
$$
(The reader is warned to remember that $u_j$ is not properly speaking a row, but two: it
is a $2\times(4k(j)-4)$ matrix and $0$ in the upper right corner of $\mathbb X_j$ is the $2\times 
2$ zero matrix.) Then \cite[Proposition 3.4]{BMS} informs us that the linearization of 
$\begin{bmatrix}
0 & P\\
P\gwia & 0
\end{bmatrix}$ is the matrix
\begin{equation}\label{eq:linearizationLP}
\mathbb L_P=\begin{bmatrix}
0 & u_1 & u_{2} & \cdots & u_{k} \\
u_{1}\gwia &\mathbb Y_1 &  &  & \\
u_{2}\gwia & & \mathbb Y_2 &  & \\
\vdots &  &  & \ddots & \\
u_{k}\gwia&  &  &  & \mathbb Y_k
\end{bmatrix}.
\end{equation}
To conclude,
\begin{equation}
(b_1-\mathbb L_P)^{-1}=\begin{bmatrix}
\left(b-\begin{bmatrix}0 & P\\
P\gwia & 0
\end{bmatrix}\right)^{-1} & \ast &  \cdots & \ast \\
\ast &\ast & \cdots & \ast \\
\vdots & \vdots & \ddots & \vdots & \\
\ast& \ast & \cdots &  \ast
\end{bmatrix}.\nonumber
\end{equation}
Here we should again remember that the only condition required is that 
$b-\begin{bmatrix}0 & P\\
P\gwia & 0
\end{bmatrix}$ is invertible, so, since $\begin{bmatrix}0 & P\\
P\gwia & 0
\end{bmatrix}$ is self-adjoint, the requirement that $\Im b>0$ will do. 

In order to apply the iteration procedure from \cite{BMS}, we only need now to 
split $\mathbb L_P$ into a sum in which elements coming from one algebra are
grouped in one matrix. Since self-adjointness is preserved by this procedure, the
subordination result of \cite{BMS} applies.

\section{Brown measure of polynomials in free variables}\label{section:polynomials}

\subsection{Identifying the Brown measure}
The linearization procedure described above guarantees that the Brown measure of
a polynomial $P$ in free variables can be expressed in terms of the $\star$-distributions
of its variables in an explicit manner. However, we would like to emphasize that the
knowledge of a significant part of the $\star$-distributions of the variables in question {\em is} 
needed: we cannot hope to obtain in general the Brown measure of $P$ out of the Brown 
measures of its variables. The knowledge of these $\star$-distributions is however guaranteed,
as noted in Section \ref{2.3}, by the knowledge of the $\G_{X\otimes 1_n}$, $n\in\mathbb N$.
As one can see easily, 
it is not in fact necessary to know $\G_{X\otimes 1_n}$ for all $n\in
\mathbb N$, but just up to a certain $n_0$ depending on the degree of 
the polynomial $P$.

An important special case is of course when all the variables $x_i$ are normal 
(for example, self-adjoint or unitary): in this case the $\star$-distribution of $x_i$ is 
determined in terms of its Brown measure (which is now nothing but
the trace applied to the spectral distribution according to the spectral
theorem). Thus our linearization procedure gives us, in particular, a way to calculate
the Brown measure of any polynomial in free self-adjoint variables out of the distribution of the variables.

Having provided the general machinery for dealing with Brown measures there are now various obvious questions to address:
\begin{itemize}
\item
Are there special cases where we can derive explicit solutions?
\item
How can we implement our algorithm to calculate numerically Brown
measures for general polynomials? Can we control 
the speed or accuracy of these calculations?
\item
Can we derive qualitative analytic features of the Brown measures?
\end{itemize}

The second question, on numerical implementation, will be addressed
somewhere else (see \cite{Spe-pol,HMS} for some preliminary results); here, we want to 
concentrate on the more analytic questions and show how we can indeed get some quite non-trivial statements out of
our general method.

As already in the self-adjoint case, there are actually not many non-trivial cases which allow an explicit description of 
Brown measures. One prominent example where one has indeed some explicit analytic formula is the case of $R$-diagonal 
elements. The Brown measure for those was calculated by Haagerup and Larsen in \cite{HaagerupLarsen}. We will show in 
the next section how their formula can be rederived in our framework.
In Section \ref{sec:triangular}. we will address another situation where an explicit calculation is possible. 
There we will consider elliptic triangular operators, which describe the limit of special Gaussian random matrix models. Only special cases of this were known before.

In this section, however, we want to start the analytic investigation of the simplest polynomial, namely the sum of two variables. 
Thus we want to address the question: what can we say about the Brown measure of $x+y$
where $x$ and $y$ are $\star$-free, given the $\star$-distribution of $x$ and the
$\star$-distribution of $y$. 

In the case where $x$ and $y$ are self-adjoint this is one
of the first fundamental questions which has been treated quite exhaustively in free probability theory, with a long list of contributions, see for example
\cite{BV-reg}. One should note that already in the case where $x=a$ and $y=ib$, with $a$ and $b$ self-adjoint 
(thus we are asking for the Brown measure of $a+ib$ where the real part $a$ and the imaginary part $b$ are free) there have been
up to now no general results on the Brown measure.

\subsection{Brown measure of the sum of two $\star$-free variables}
The Hermitization and linearization method from the last section show that
the treatment of an arbitrary polynomial in two $\star$-free variables
requires the analysis of the matrix $\mathbb L_P$ from \eqref{eq:linearizationLP}, split into the sum of the two terms corresponding to the two free variables.

At the moment we are not able to analyze the analytic features of this general framework; what we can and will 
do here is to treat in some detail the case where $k=1$, i.e., the case of a linear polynomial, corresponding to the Brown measure of the sum 
of two $\star$-free random variables $x$ and $y$.
Note that in this situation there is no need for a linearization and we just have
to understand the $\mathcal M_2(\mathbb C)$-valued distribution of
$$\begin{bmatrix} 0 & x+y \\ (x+y)\gwia & 0 \end{bmatrix}=X+Y=
\begin{bmatrix} 0 & x \\ x\gwia & 0 \end{bmatrix}+
\begin{bmatrix} 0 & y \\ y\gwia & 0 \end{bmatrix}$$
in terms of the $\mathcal M_2(\mathbb C)$-valued distributions of $X$ and of $Y$. Note that $X$ and $Y$ are free over $\mathcal M_2(\mathbb C)$

Let us first remind the reader of the result from 
\cite{coalgebra,BMS} related to subordination: there exist two analytic self-maps
of the upper half-plane of $\mathcal M_2(\mathbb C)$, called $\omega_1,\omega_2$,
so that
\begin{equation}\label{subord}
(\omega_1(b)+\omega_2(b)-b)^{-1}=\G_X(\omega_1(b))=\G_Y(\omega_2(b))=\G_{X+Y}(b),
\end{equation}
for all $b\in\mathcal M_2(\mathbb C)$ with $\Im b>0$. We shall be concerned with a special
type of $b$, namely $b=\begin{bmatrix}z & \lambda\\
\overline{\lambda} & w
\end{bmatrix}$. As mentioned before, while the correspondences in $z$ and $w$ are analytic, the 
correspondence in $\lambda$ is not, so we shall view it as a parameter, on which however
there is a continuous corresponcence when $z$ and $w$ in the upper half plane of $\mathbb C$ are fixed. Thus, we denote below
$\omega_j(b), \G(b),\F(b)$ by $\omega_j(z,w), \G(z,w),\F(z,w)$, sometimes
adding the parameter $\lambda$, when it becomes important in terms of our analysis (it will
usually be fixed).

For a fixed $\lambda\in\mathbb C$, we are interested in the behavior of $\G_{X+Y}(z,w)$
close to $z=w=0$. For our purposes, it will be enough to consider the case $z=w$ and view the 
functions involved as single-variable holomorphic functions. For obvious practical purposes,
we would like to argue that $\lim_{\epsilon\to0}\G_{X+Y}(i\epsilon,i\epsilon)$ exists for all 
$\lambda \in\mathbb C$. Sadly, this is, quite trivially, not true. We shall analyze this problem in 
several steps, obtaining along the way side results which we believe interesting in their own right. 

\subsubsection{The question of left and right invariant projections}\label{sec:proj}
We remind the reader of some general facts about elements in finite von Neumann algebras (i.e., when we have a faithful normal trace) (see \cite{Stratila-Zsido}). 
\begin{itemize}

\item
If $a\in\mathfrak A$ is an arbitrary, non-Hermitian, element,
and $p=p^*=p^2\in\mathfrak A\setminus\{0\}$ is a projection so that $ap=\lambda p$, then 
$p\leq\ker((a-\lambda)\gwia(a-\lambda))$ (we denote here by $\ker(b)$ both the (closed) space on which 
$b$ is zero and the orthogonal projection onto this space). Indeed, recalling that $\phi$ is
a faithful trace and $p(a-\lambda)\gwia(a-\lambda)p\ge0$,
$$
0=\phi\left((a-\lambda)\gwia(a-\lambda)p\right)=\phi\left(p(a-\lambda)\gwia(a-
\lambda)p\right),
$$
 implies that $\ker((a-\lambda)\gwia(a-\lambda))\ge p$. Similarly, if $p\leq
\ker((a-\lambda)\gwia(a-\lambda))$, then by the above equality $(a-\lambda)p=0$.

\item If $a=v(a^*a)^\frac12=v|a|$ is the polar decomposition of $a$, then the partial
isometry $v$ can be completed to a unitary operator in the von Neumann algebra generated 
by $a$.

\item If $p=\ker(a-\lambda)$, then there exists $q=q\gwia=q^2\in\mathfrak A$ such that 
$\phi(p)=\phi(q)$ and $qa={\lambda}q$. Indeed, by replacing $a$ with $a-\lambda$, we
may assume that $\lambda=0$. We assume that the partial isometry in the polar decompostion of
$a$ is (completed to) a unitary operator. Then 
$$
ap=0\implies |a|p=0\text{ and }p|a|=0.
$$
In particular, if we let $q=vpv\gwia$, then $qa=vpv\gwia a=vpv\gwia v|a|=vp|a|=0$. Since $q^2=vp
v\gwia vpv\gwia=vp^2v\gwia=vpv\gwia=q$ and $\phi(q)=\phi(vpv\gwia)=\phi(v\gwia vp)=\phi(p)$, we 
conclude that $q$ is a projection equivalent to $p$.

\item Thus, if $p=\ker(a-\lambda)$, $a=v|a|$, $q=vpv\gwia$, then
\begin{equation}
ap=\lambda p,\quad qa=\lambda q,\quad pa\gwia=\overline{\lambda}p,\quad a\gwia
q=\overline{\lambda}q.
\end{equation}

\end{itemize}
Now we are in the position to carry on almost to the letter the analysis from \cite{BV-reg}. 
We shall denote $\G_{X+Y}(i\epsilon,i\epsilon)$ simply by $\G_{X+Y}(i\epsilon)$, and 
similarly for the other functions involved. We have
$$
\lim_{\epsilon\downarrow0}\G_{X+Y}(i\epsilon)=\begin{bmatrix}
\displaystyle\lim_{\epsilon\downarrow0} g_{X+Y,11}(i\epsilon) &
\displaystyle\lim_{\epsilon\downarrow0} g_{X+Y,12}(i\epsilon)\\
\displaystyle\lim_{\epsilon\downarrow0} g_{X+Y,21}(i\epsilon) &
\displaystyle\lim_{\epsilon\downarrow0} g_{X+Y,22}(i\epsilon)
\end{bmatrix},
$$
where we remind the reader that 
\begin{eqnarray*}
g_{X+Y,11}(i\epsilon) & = & -i\epsilon\phi\left(\big((\lambda-x-y)(\lambda-x-y)\gwia+\epsilon^2
\big)^{-1}\right)\\
g_{X+Y,12}(i\epsilon) &=& \phi\left((\lambda-x-y)\big((\lambda-x-y)\gwia (\lambda-x-y)+\epsilon^2
\big)^{-1}\right) \\
g_{X+Y,21}(i\epsilon)&=&\phi\left((\lambda-x-y)\gwia\big((\lambda-x-y)(\lambda-x-y)\gwia+
\epsilon^2\big)^{-1}\right)\\
g_{X+Y,22}(i\epsilon) & = & -i\epsilon\phi\left(\big( (\lambda-x-y)\gwia(\lambda-x-y)+\epsilon^2 
\big)^{-1}\right).
\end{eqnarray*}
Since $\phi$ is a trace, $g_{X+Y,11}(i\epsilon)=g_{X+Y,22}(i\epsilon)\in-i(0,+\infty)$ and 
$g_{X+Y,12}(i\epsilon)=\overline{g_{X+Y,21}(i\epsilon)}$. Clearly, the
same holds if $X+Y$ is replaced by $X$ or $Y$.
We assume that $0<p=\ker(\lambda-x-y)<1$, where $x,y\in\mathfrak A\setminus\mathbb C$ are
$\star$-free. Recall that the last hypothesis makes $\mathcal M_2(\mathbb C\left\langle x,x\gwia
\right\rangle)$ and $\mathcal M_2(\mathbb C\left\langle y,y\gwia
\right\rangle)$ free over $\mathcal M_2(\mathbb C)$.
For convenience, we shall temporarily denote $a=\lambda-x-y$, so that $ap=a\gwia ap=0$
and 
$$g_{X+Y,11}(i\epsilon)=-i\epsilon\phi\left((aa\gwia+\epsilon^2)^{-1}\right),\quad
g_{X+Y,12}(i\epsilon)=\phi\left(a(a\gwia a+\epsilon^2)^{-1}\right),$$ 
$$g_{X+Y,22}(i\epsilon)=-i\epsilon\phi\left((a\gwia a+\epsilon^2)^{-1}\right),\quad
g_{X+Y,21}(i\epsilon)=\phi\left(a\gwia(aa\gwia+\epsilon^2)^{-1}\right).$$ 
The weak limits
\begin{eqnarray*}
\lim_{\epsilon\downarrow0}i\epsilon\left(-i\epsilon(a\gwia a+\epsilon^2)^{-1}\right) & = &
\lim_{\epsilon\downarrow0}\epsilon^2(a\gwia a+\epsilon^2)^{-1}\\
& = & 1-\lim_{\epsilon\downarrow0}a\gwia a(a\gwia a+\epsilon^2)^{-1}\\
& = & \ker(a\gwia a)\\
&=&p,
\end{eqnarray*}
and 
$$
\lim_{\epsilon\downarrow0}i\epsilon\left(-i\epsilon(aa\gwia+\epsilon^2)^{-1}\right)=
\lim_{\epsilon\downarrow0}\epsilon^2(aa\gwia+\epsilon^2)^{-1}=\ker(aa\gwia)=q
$$
hold, so that in particular 
$$
\lim_{\epsilon\downarrow0}i\epsilon g_{X+Y,11}(i\epsilon)=\lim_{\epsilon\downarrow0}i\epsilon 
g_{X+Y,22}(i\epsilon)=\phi(p)=\phi(q).
$$
For the $(1,2)$ and $(2,1)$ entries, the situation is slightly more delicate: for the polar decomposition 
$a=v (a\gwia a)^\frac12$, 
using the generalized Schwarz
inequality $|\phi(x\gwia y)|^2\leq\phi(x\gwia x)\phi(y\gwia y)$
applied to $x=v$
and $y=(a\gwia a)^\frac12\big(a\gwia a+\epsilon^2\big)^{-1}$, we write
$$
\left|\phi\left(a\big(a\gwia a+\epsilon^2\big)^{-1}\right)\right|\leq\left[\phi(1)\right]^\frac{1}2
\left[\phi\left(a\gwia a\big(a\gwia a+\epsilon^2\big)^{-2}\right)\right]^\frac12.
$$
Since $2\epsilon^2t(t+\epsilon^2)^{-2}<1$ for all $t\ge0$, we obtain by applying continuous 
functional calculus to $\epsilon^2a\gwia a\big(a\gwia a+\epsilon^2\big)^{-2}$ and by the positivity of
$\phi$ that 
$$
\left|\epsilon\phi\left(a\big(a\gwia a+\epsilon^2\big)^{-1}\right)\right|\leq\left[\phi(1)\right]^\frac{1}2
\left[\phi\left(\epsilon^2a\gwia a\big(a\gwia a+\epsilon^2\big)^{-2}\right)\right]^\frac12<1.
$$
On the other hand, observe that $\lim_{\epsilon\downarrow0}
\epsilon^2t(t+\epsilon^2)^{-2}=0$ pointwise for $t\in[0,+\infty)$, so, if $\theta$ denotes the
distribution of $a\gwia a$ with respect to $\phi$, then by dominated convergence we obtain
$$
\lim_{\epsilon\downarrow0}\phi\left(\epsilon^2a\gwia a\big(a\gwia a+\epsilon^2\big)^{-2}\right)
=\lim_{\epsilon\downarrow0}\int_{[0,+\infty)}\epsilon^2t(t+\epsilon^2)^{-2}\,d\theta(t)=0.
$$
We conclude that, as expected,
$$
\lim_{\epsilon\downarrow0}i\epsilon g_{X+Y,21}(i\epsilon)=\lim_{\epsilon\downarrow0}i\epsilon
g_{X+Y,12}(i\epsilon)=0.
$$
Recall from \eqref{subord} that 
\begin{equation}\label{subord-prime}
\omega_1(i\epsilon)+\omega_2(i\epsilon)=\begin{bmatrix}
i\epsilon&\lambda\\
\overline{\lambda}&i\epsilon
\end{bmatrix}+\F_{X+Y}(i\epsilon).
\end{equation}
The function $\F_{X+Y}$ is easily obtained by inverting $\G_{X+Y}$:
$$
\F_{X+Y}=\frac{1}{g_{X+Y,11}g_{X+Y,22}-g_{X+Y,12}g_{X+Y,21}}\begin{bmatrix}
g_{X+Y,22} & -g_{X+Y,12} \\
-g_{X+Y,21} & g_{X+Y,11}
\end{bmatrix}.
$$
It follows from the corresponding property of the matrix-valued Cauchy transform $\G$ that $f_{X+Y,11}(i\epsilon)=f_{X+Y,22}(i\epsilon)\in i(0,+\infty)$, and $f_{X+Y,12}(i\epsilon)
=\overline{f_{X+Y,21}(i\epsilon)}$, whith the same holding when $X+Y$ is replaced by $X$ or $Y$.
Consider the entrywise limits $\lim_{\epsilon\downarrow0}\F_{X+Y}(i\epsilon)/i\epsilon$ in the above.
Using the previously obtained estimates,
\begin{eqnarray*}
\lefteqn{\lim_{\epsilon\downarrow0}\frac{f_{X+Y,11}(i\epsilon)}{i\epsilon} }\\
& = &\lim_{\epsilon\downarrow0}\frac{[i\epsilon g_{X+Y,22}(i\epsilon)]}{[i\epsilon g_{X+Y,11}(i\epsilon)]
[i\epsilon g_{X+Y,22}(i\epsilon)]-[i\epsilon g_{X+Y,12}(i\epsilon)][i\epsilon g_{X+Y,21}(i\epsilon)]}\\
& = & \frac{\phi(p)}{\phi(q)\phi(p)-0\cdot0}=\frac{1}{\phi(q)},
\end{eqnarray*}
and
\begin{eqnarray*}
\lefteqn{\lim_{\epsilon\downarrow0}\frac{f_{X+Y,12}(i\epsilon)}{i\epsilon} }\\
& = &\lim_{\epsilon\downarrow0}\frac{-[i\epsilon g_{X+Y,12}(i\epsilon)]}{[i\epsilon g_{X+Y,11}(i\epsilon)]
[i\epsilon g_{X+Y,22}(i\epsilon)]-[i\epsilon g_{X+Y,12}(i\epsilon)][i\epsilon g_{X+Y,21}(i\epsilon)]}\\
& = & \frac{-0}{\phi(q)\phi(p)-0\cdot0}=0.
\end{eqnarray*}
Identical computations provide $\lim_{\epsilon\downarrow0}f_{X+Y,11}(i\epsilon)/i\epsilon=1/\phi(p)=
1/\phi(q)$ and $\lim_{\epsilon\downarrow0}f_{X+Y,21}(i\epsilon)/i\epsilon=0$.
In $\mathcal M_2(\mathbb C)$-norm, we obtain all of
$$
\lim_{\epsilon\downarrow0}\F_{X+Y}(i\epsilon)=0,\quad\lim_{\epsilon\downarrow0}
\frac{1}{i\epsilon}\F_{X+Y}(i\epsilon)=\lim_{\epsilon\downarrow0}\frac1\epsilon\Im\F_{X+Y}(i\epsilon)=
\phi(p)^{-1}1_2.
$$
Recall that $\omega_1,\omega_2$ map the upper half-plane of $\mathcal M_2(\mathbb C)$ into itself,
increase the imaginary part (meaning that $\Im\omega_k(b)\ge\Im b$ whenever $b\in\mathcal M_2(\mathbb C),\Im b>0$),
and $\Im\omega_k(b)\le\Im\F_{X+Y}(b)$ if $\Im b>0$, $k=1,2$ - see \cite{coalgebra,BMS,BPV}). This,
together with the behavior of $\F_{X+Y}$ near zero, allows us to conclude that 
\begin{equation}\label{Jim}
\lim_{\epsilon\downarrow0}\frac1\epsilon\Im\omega_1(i\epsilon),
\lim_{\epsilon\downarrow0}\frac1\epsilon\Im\omega_2(i\epsilon)\in
\left[1_2,\frac{1}{\phi(p)}1_2\right]. %%%%%%Removed the 1+, as it was clearly superfluous
\end{equation}
 and  Ky Fan's operator generalization of the Julia-Carath\'eodory 
Theorem from \cite{Fan} allows us to conclude that $\lim_{\epsilon\downarrow0}\omega_j(i\epsilon)$, $j=1,2$,
exist and are selfadjoint $2\times2$ matrices. Recall that both $\F_X$ and $\F_Y$ preserve the set of
matrices of the form  $\begin{bmatrix} i\epsilon & \lambda \\ \overline{\lambda} & i\epsilon\end{bmatrix}$, $\epsilon>0,\lambda\in\mathbb C$.
Then, according to \cite[Theorem 2.7]{BMS}, so do $\omega_1$ and $\omega_2$.
We conclude that
\begin{equation}\label{ady}
\lim_{\epsilon\downarrow0}\omega_1(i\epsilon)=\begin{bmatrix}
0 & u_1\\
\overline{u}_1&0\end{bmatrix},\quad\lim_{\epsilon\downarrow0}\omega_2(i\epsilon)=\begin{bmatrix}
0 & u_2\\
\overline{u}_2&0\end{bmatrix},
\end{equation}
where $u_1+u_2=\lambda$, and
\begin{equation}\label{corner}
w_k:=\lim_{\epsilon\downarrow0}\frac{\Im\omega_k(i\epsilon)}{\epsilon}=
\lim_{\epsilon\downarrow0}\frac{1}{i\epsilon}\left(\omega_k(i\epsilon)-\begin{bmatrix}
0 & u_k\\
\overline{u}_k&0\end{bmatrix}\right)\in\left[1_2,\frac{1}{\phi(p)}1_2\right],
\end{equation}
for $ k\in\{1,2\},$ all being norm limits.
Subtracting $\begin{bmatrix}0&\lambda\\ \overline{\lambda} & 0\end{bmatrix}$
from both sides of \eqref{subord-prime} then dividing by $i\epsilon$ and then letting $\epsilon$
tend to zero guarantees that 
\begin{equation}\label{Cara}
\lim_{\epsilon\downarrow0}\frac{(\omega_1(i\epsilon)+\omega_2(i\epsilon))_{12}-\lambda}{i\epsilon}
=0,
\end{equation}
with a similar result for the $(2,1)$ entry, except that $\lambda$ is replaced by 
$\overline{\lambda}$.

We use next the
coalgebra morphism property of the conditional expectation proved by Voiculescu in \cite{coalgebra}: it 
is known that whenever $X,Y$ are free over $B$, one has 
$\mathbb E_{B\langle X\rangle}\left[(b-X-Y)^{-1}\right]
=(\omega_1(b)-X)^{-1}$. We apply this to our $2\times 2$ matrices $X,Y$ and $B=\mathcal M_2(\mathbb 
C)$ to write
\begin{eqnarray*} \lefteqn{
\mathbb E_{\mathcal M_2(\mathbb C)\langle X\rangle}\left[\left(\begin{bmatrix} i\epsilon & \lambda\\ \overline{\lambda} & i\epsilon\end{bmatrix}-X-Y\right)^{-1}\right]}\quad\quad\quad\quad
\quad\quad\quad\quad\quad\quad\\
&=&\begin{bmatrix}
(\omega_1(i\epsilon))_{11} & (\omega_1(i\epsilon))_{12}-x\\
(\omega_1(i\epsilon))_{21}-x\gwia & (\omega_1(i\epsilon))_{22}
\end{bmatrix}^{-1}.
\end{eqnarray*}
The  formulas for the (1,1) and (1,2) entries of the right-hand matrix are 
$$
(\omega_1(i\epsilon))_{22}
\left[(\omega_1(i\epsilon))_{11}(\omega_1(i\epsilon))_{22}-((\omega_1(i\epsilon))_{12}-x)((\omega_1(i\epsilon))_{21}-x\gwia)\right]^{-1}
$$
and
$$\left[{((\omega_1(i\epsilon))_{12}-x)((\omega_1(i\epsilon))_{21}-x\gwia)}-(\omega_1(i\epsilon))_{11}{(\omega_1(i\epsilon))_{22}}\right]^{-1}((\omega_1(i\epsilon))_{12}-x).
$$
The assumption
$\ker(x+y-\lambda)= p$ implies that 
$$
\left(X+Y-\begin{bmatrix} 0 &\lambda \\ \overline{\lambda} & 0 \end{bmatrix}\right)\begin{bmatrix} q &0\\ 0 & p \end{bmatrix}= \begin{bmatrix} q &0\\ 0 & p \end{bmatrix}\left(X+Y-\begin{bmatrix} 0 &\lambda \\ \overline{\lambda} & 0 \end{bmatrix}\right)=0.
$$
 Thus, the weak limit
$$
\lim_{\epsilon\downarrow0}\begin{bmatrix} i\epsilon &0 \\ 0 & i\epsilon \end{bmatrix}
\left(\begin{bmatrix} i\epsilon &\lambda \\ \overline{\lambda} & i\epsilon\end{bmatrix}-X-Y\right)^{-1}=
\begin{bmatrix} q &0\\ 0 & p \end{bmatrix}
$$
holds. Recalling the weak continuity of the conditional expectation (w.r.t. $\phi$) onto a von Neumann 
subalgebra of $\mathfrak A$, we obtain 
\begin{eqnarray}
\lefteqn{\lim_{\epsilon\downarrow0}\begin{bmatrix} i\epsilon &0 \\ 0 & i\epsilon \end{bmatrix}
\begin{bmatrix}
(\omega_1(i\epsilon))_{11} & (\omega_1(i\epsilon))_{12}-x\\
(\omega_1(i\epsilon))_{21}-x\gwia & (\omega_1(i\epsilon))_{22}\end{bmatrix}^{-1}}\nonumber
\quad\quad\quad\quad\quad\quad\quad\quad\quad\quad\quad\quad\quad\quad\quad\quad\\
& = &
\begin{bmatrix} \mathbb E_{\mathbb C\langle x,x\gwia\rangle}[q] &0\\ 0 & 
\mathbb E_{\mathbb C\langle x,x\gwia\rangle}[p] \end{bmatrix}.\label{proiectie}
\end{eqnarray}
A similar relation is deduced for $\omega_2$ and $y$.
Using the formulas for the entries of the term under the limit, we obtain
\begin{eqnarray*}
\lefteqn{\mathbb E_{\mathbb C\langle x,x\gwia\rangle}[q]=\lim_{\epsilon\downarrow0}i\epsilon
(\omega_1(i\epsilon))_{22}}\\
& &\mbox{}\times
\left[(\omega_1(i\epsilon))_{11}(\omega_1(i\epsilon))_{22}-((\omega_1(i\epsilon))_{12}-x)((\omega_1(i\epsilon))_{21}-x\gwia)\right]^{-1}\\
& = & \lim_{\epsilon\downarrow0}\frac{i\epsilon}{(\omega_1(i\epsilon))_{11}}\cdot{(\omega_1(i\epsilon))_{22}}{(\omega_1(i\epsilon))_{11}}\\
& &\mbox{}\times
\left[(\omega_1(i\epsilon))_{11}(\omega_1(i\epsilon))_{22}-((\omega_1(i\epsilon))_{12}-x)((\omega_1(i\epsilon))_{21}-x\gwia)\right]^{-1}\\
& = & (w_1)_{11}^{-1}\cdot\lim_{\epsilon\downarrow0}{(\omega_1(i\epsilon))_{22}}{(\omega_1(i\epsilon))_{11}}\\
& &\mbox{}\times
\left[(\omega_1(i\epsilon))_{11}(\omega_1(i\epsilon))_{22}-((\omega_1(i\epsilon))_{12}-x)((\omega_1(i\epsilon))_{21}-x\gwia)\right]^{-1}.
\end{eqnarray*}
 Thus,
\begin{eqnarray}
\lefteqn{(w_1)_{11}\mathbb E_{\mathbb C\langle x,x\gwia\rangle}[q]=-
\lim_{\epsilon\downarrow0}\left[(\omega_1(i\epsilon))_{11}(\omega_1(i\epsilon))_{22}\times\frac{}{}
\right.}\nonumber\\
& & \left.\left[((\omega_1(i\epsilon))_{12}-x)((\omega_1(i\epsilon))_{21}-x\gwia)-
(\omega_1(i\epsilon))_{11}(\omega_1(i\epsilon))_{22}\right]^{-1}\frac{}{}\right]\label{29}
\end{eqnarray}
The operator $\mathbb E_{\mathbb C\langle x,x\gwia\rangle}[q]$ is nonzero, nonnegative and bounded 
from above by $1$. Similarly, considering the (1,2) entry, we have the weak limit
\begin{eqnarray*}
\lefteqn{\lim_{\epsilon\downarrow0}i\epsilon\left[{((\omega_1(i\epsilon))_{12}-x)((\omega_1
(i\epsilon))_{21}-x\gwia)}-(\omega_1(i\epsilon))_{11}{(\omega_1(i\epsilon))_{22}}\right]^{-1}}\\
& & \mbox{}\times%(\omega_1(i\epsilon))_{22} Removed it, it's a typo, see previous page
((\omega_1(i\epsilon))_{12}-x)=0.
\quad\quad\quad\quad\quad\quad\quad\quad\quad\quad\quad\quad\quad
\end{eqnarray*}
Observe that 
\begin{eqnarray*}
\lefteqn{
\|((\omega_1(i\epsilon))_{12}-x)((\omega_1(i\epsilon))_{21}-x\gwia)-(u_1-x)(u_1-x)\gwia\|}\quad\quad
\quad\quad\quad\quad\quad\quad\\
&\leq&((\omega_1(i\epsilon))_{12}+|u_1|+2\|x\|)|(\omega_1(i\epsilon))_{12}-u_1|.%%%%Changed here from 2(1+\|x\|)|(\omega_1(i\epsilon))_{12}-u_1| to the precise estimate
\end{eqnarray*}
Equations \eqref{Jim} and \eqref{corner} guarantee that $\lim_{\epsilon\downarrow0}
|(\omega_1(i\epsilon))_{12}-u_1|/\epsilon<1/\phi(p).$ %%%%%Cut the 1+

As $(\omega_1(i\epsilon))_{12}$ converges to $u_1$ and $x$ is constant, the above allows us to 
conclude that $\mathbb E_{\mathbb C\langle x,x\gwia\rangle}[q](u_1-x)=0.$ A similar
computation, in which we use entries (2,2) and (2,1) of the corresponding matrices, provides 
$({u}_1-x)\mathbb E_{\mathbb C\langle x,x\gwia\rangle}[p]=0$. Considering the support projection 
of $\mathbb E_{\mathbb C\langle x,x\gwia\rangle}[q]$ and $p_1$ of 
$\mathbb E_{\mathbb C\langle x,x\gwia\rangle}[p]$, we find that $x-u_1$ has nonzero left and right 
invariant projections (with $p_1$ being the right-invariant projection). A similar argument shows that 
$y-u_2$ has left and right invariant projections, with $p_2$ being the right-invariant projection.

The (sum of the) length of these projections is deduced the following way: it is clear that $p\geq
p_1\wedge p_2.$ 
Recalling the 
subordination relation and \eqref{proiectie}, we obtain 
$$
\phi(p)^{-1}=%(\omega_1'(0))
(w_1)_{11}+(%\omega_2'(0)
w_2)_{11}-1
$$ 
(recall that $\omega_k(i\epsilon)$ have equal entries on the
diagonal, so picking $p$ and $(1,1)$ or $q$ and (2,2) makes no difference). On the other hand, we
recall that (by definition), $\mathbb E_{\mathbb C\langle x,x\gwia\rangle}$ preserves the trace:
$\phi(p)=\phi(\mathbb E_{\mathbb C\langle x,x\gwia\rangle}[p])$. Since $0\leq
\mathbb E_{\mathbb C\langle x,x\gwia\rangle}[p]\leq1$, it follows immediately from elementary
functional calculus that the trace of the support of $\mathbb E_{\mathbb C\langle x,x\gwia\rangle}[p]$
cannot be less than $\phi(\mathbb E_{\mathbb C\langle x,x\gwia\rangle}[p])$, so $\phi(p_1),\phi(p_2)
\ge \phi(p)$. Also, by applying $\phi$ in relation \eqref{29} it follows that 
$$\phi(p)(\omega_1'(0))_{11}=
\phi(\mathbb E_{\mathbb C\langle x,x\gwia\rangle}[p])(%\omega_1'(0)
w_1)_{11}\leq\phi(p_1).$$ Clearly,
a similar relation holds for $p_2$, $y$ and $(%\omega_2'(0)
w_2)_{11}$. Thus, 
$$
\frac{1}{\phi(p)}=(%\omega_1'(0)
w_1)_{11}+(%\omega_2'(0)
w_2)_{11}-1\leq\frac{\phi(p_1)+\phi(p_2)}{\phi(p)}-1,
$$
which is equivalent to
$$
\phi(p_1)+\phi(p_2)\ge\phi(p)+1.
$$
This, together with the relation $p\geq p_1\wedge p_2$, implies $\phi(p_1)+\phi(p_2)=\phi(p)+1$.

Thus we have proved the following result, paralleling \cite[Theorem 7.4]{BV-reg}.
\begin{proposition}
If $x,y\in\mathfrak A\setminus\mathbb C1$ are $\star$-free with respect to $\phi$ and there exist a 
projection $p\in\mathfrak A\setminus\{0\}$, $\lambda\in\mathbb C$, so that $(x+y)p=\lambda p$, then
\begin{enumerate}
\item $p=\ker((x+y-\lambda)\gwia(x+y-\lambda))$;
\item there exist $p_1,p_2$ projections in $\mathfrak A$ and $u_1,u_2\in\mathbb C$
so that 
\begin{itemize}
\item $xp_1=u_1p_1$ and $yp_2=u_2p_2$.
\item $u_1+u_2=\lambda$.
\item $\phi(p_1)+\phi(p_2)=\phi(p)+1.$
\end{itemize}
\end{enumerate}
Conversely, if the three conditions of item $(2)$ above hold, then $p:=p_1\wedge p_2$ satisfies
$(x+y)p=(u_1+u_2)p$.
\end{proposition}

\begin{remark}
In \cite{BV-reg} it is shown that under the same hypotheses, if $x=x\gwia$ and $y=y\gwia$,
then $\omega_1'(\lambda)\mathbb E_{\mathbb C\langle x\rangle}[p]$ is a projection.
It would be interesting to determine whether $\mathbb E_{\mathbb C\langle x,x\gwia\rangle}[q]$ and $
\mathbb E_{\mathbb C\langle x,x\gwia\rangle}[p]$ are themselves multiples of projections.
\end{remark}

\begin{remark}
Let us note that, regrettably, the limits $\lim_{\epsilon\downarrow0}\G_{X+Y}(i\epsilon)$ in general will 
not provide the value of the Cauchy-Stieltjes transform in $\lambda$. Indeed, assume that $x=x\gwia$ 
and $y=y\gwia$, neither a multiple of the identity. Then it is known that, roughly speaking,
$G_{x+y}$ extends continuously to the real line and $G_{x+y}(r)\in\mathbb C^-$ (the lower half plane of $\C$) has an analytic 
extension around $r$ for most points in the spectrum of $x+y$. However, for $\lambda=r$ being 
one of those points in the spectrum of $x+y$, we have
\begin{eqnarray*}
\lim_{\epsilon\downarrow0}\G_{X+Y,21}(i\epsilon)&=&
\lim_{\epsilon\downarrow0}\phi((r-x-x)((r-x-y)^2+\epsilon^2)^{-1})\\
&=&
\lim_{\epsilon\downarrow0}\int_\mathbb R\frac{r-t}{(r-t)^2+\epsilon^2}\,d\mu_{x+y}(t),
\end{eqnarray*}
which is simply the Hilbert transform of $\mu_{x+y}$ evaluated at $r$, a real number. On the
bright side, note that this limit {\em does} exist for all $r\in\mathbb R$, with the exception
of those points $r$ where $G_{x+y}$ is infinite (and it is known that there are only finitely many
such points). On the even brighter side,
the $(1,1)$ entry will provide through the same argument the imaginary part of $G_{x+y}(r)$.
Of course, it is natural that it should be so, because otherwise we would have 
$\frac{\partial}{\partial\overline{\lambda}}G_{x+y}(r)=0$ for most points $r$ in the spectrum
of $x+y$. This simple example shows us that $\lambda\mapsto\G_{X+Y}(0)$ has little
chance of being generally continuous, and, in particular, that \cite[Theorem 3.3 (3)]{B-PTRF} does 
not hold in the operator-valued context. 
\end{remark}
We dare nevertheless to make the
following conjecture.
\begin{conjecture}
If $x,y\in\mathfrak A$ are $\star$-free w.r.t. the trace state $\phi$ and the spectrum of each contains
more than one point, then
the function $\mathbb C\ni\lambda\mapsto\G_{X+Y}(0)$ is continuous when restricted to each 
component of the spectrum of $x+y$. The points of discontinuity of $\mathbb C\ni\lambda
\mapsto\G_{X+Y}(0)$ 
belong to the closure of the resolvent of $x+y$.
\end{conjecture}

\section{Brown measure of $R$--diagonal elements}
\label{sec:rdiagonal}

In this section we will show that $R$-diagonal operators fit very nicely
in the Hermitization and subordination frame and that one can recover from
this point of view in a quite systematic way the result of Haagerup and Larsen
on the Brown measure of $R$-diagonal operators.

Let us first recall, for later use, the definition of the $S$-transform.
Recall that for a random variable $y$ with
$\phi(y)\neq 0$ we put \cite{VDN}
$$ \psi_y (\lambda): = \phi\left( (1-\lambda y)^{-1} \right) -1; $$
one can show that
$$ S_y (\lambda):= \frac{\lambda+1}{\lambda} \psi^{\langle -1 \rangle}_y(\lambda) $$
is well--defined in some neighborhood of $0$; it is called
Voiculescu's $S$--transform of $y$. It was proved in
\cite{LarsenPhD,HaagerupLarsen} that if $y$ is a positive operator
then its $S$--transform has an analytic continuation to an
interval $(\mu_y\{0\}-1, 0]$; it has a strictly negative
derivative on this interval and
$$ S_y(\mu_y\{0\}-1)=\phi\left( y^{-2} \right),\qquad S_y(0)= \frac{1}{\phi(y^2)}.$$

Let now $x=ua\in(\mathfrak A,\phi)$ be $R$-diagonal \cite{NSbook}. Recall that this means that $u$ and $a$ are $\star$-free
with respect to $\phi$, and that $u$ is a Haar unitary and $a\ge0$. 
Note that $R$-diagonal operators are, in $\star$-moments, the limits of an important class of random matrices, namely of bi-unitarily invariant random matrices.

We shall denote by $\mathcal D$ the
commutative $C\gwia$-algebra of diagonal matrices in $\mathcal M_2(\mathbb C)$, where we
consider the same inclusion of $\mathcal M_2(\mathbb C)$ in $\mathcal M_2(\mathfrak A)$
as in the previous sections.
Let us recall from \cite{NSS} that an element $x$  in a tracial $W\gwia$-noncommutative
probability space is $R$-diagonal if and only if the matrix \begin{tiny}$\begin{bmatrix}0&x\\
x\gwia&0\end{bmatrix}$\end{tiny} is free over
$\mathcal D$ from $\mathcal M_2(\mathbb C)$ with respect to the expectation
$$\mathbb E_\mathcal D\colon\mathcal M_2(\mathfrak A)\to
\mathcal D,\qquad
\mathbb E_\mathcal D\begin{bmatrix}a_{11}&a_{12}\\
a_{21}&a_{22}\end{bmatrix}=\begin{bmatrix}\phi(a_{11})&0\\
0&\phi(a_{22})\end{bmatrix}.$$ 
Fix now a $\lambda$ in the upper half plane of $\mathbb C$. With the notation $\mathcal D\left\langle
X\right\rangle$ for the $\star$-algebra generated by $\mathcal D$ and $X$, it is
obvious that 
$$\mathcal D\left\langle\begin{bmatrix}0&\lambda\\
\overline{\lambda}&0\end{bmatrix}\right\rangle=\mathcal M_2(\mathbb C).$$ 
(The reader can
easily verify this by, for example, constructing all matrix units out of the two nontrivial
projections of $\mathcal D$ and the element 
\begin{tiny}$\begin{bmatrix}0&\lambda\\
\overline{\lambda}&0\end{bmatrix}$.\end{tiny})
The fundamental result \cite[Theorem 3.8]{coalgebra} of Voiculescu is written in this context as:
\begin{equation}\label{R-dia}
\mathbb E
\left[\left(\begin{bmatrix}z&\lambda\\
\overline{\lambda}&w\end{bmatrix}-\begin{bmatrix}0&x\\
x\gwia&0\end{bmatrix}\right)^{-1}\right]=\begin{bmatrix}\omega_1(z,w)&\lambda\\
\overline{\lambda}&\omega_2(z,w)\end{bmatrix}^{-1},
\end{equation}
where 
$$\mathbb E=
\mathbb E_{\tiny{\mathcal D\left\langle\begin{bmatrix}0&\lambda\\
\overline{\lambda}&0\end{bmatrix}\right\rangle}}=\mathbb E_{\mathcal M_2(\mathbb C)}$$ 
denotes, as above, the unique conditional
expectation onto $\mathcal M_2(\mathbb C)$ which preserves the trace, and which is given by evaluation of $\phi$ on the entries of the matrix. As before,
the $(2,1)$ entry is the object of interest $G_{\mu_x}$, determined by the functions $\omega_1$,
$\omega_2$ via a straightforward algebraic relation. (It is obvious that the functions $\omega_1$
and $\omega_2$ depend also on $\lambda$, and this dependence is relevant to us. While it
might be unfair to the reader, we will follow tradition and suppress this dependence in notations.)

On the other hand, the subordination function $\omega=\begin{bmatrix}\omega_1&0\\
0&\omega_2\end{bmatrix}$, obtained by applying $\mathbb E_\mathcal D$ in \eqref{R-dia}, is
determined by \cite[Theorem 2.2]{BMS} via the iteration procedure, if desired, or via
straightforward direct computation in terms of (functions derived from) the Cauchy-Stieltjes
transform of $a^2$. Let $w=z$.
Performing the inversion in the right hand side of \eqref{R-dia} gives
\begin{multline}
\label{trente-neuf}
\begin{bmatrix}\omega_1(z,z)&\lambda\\
\overline{\lambda}&\omega_2(z,z)\end{bmatrix}^{-1}=\\
\begin{bmatrix}\omega_2(z,z)\left(\omega_1\omega_2(z,z)-|\lambda|^2\right)^{-1}&\lambda\left(|\lambda|^2-\omega_1\omega_2(z,z)\right)^{-1}\\
\overline{\lambda}\left(|\lambda|^2-\omega_1\omega_2(z,z)\right)^{-1}&\omega_1(z,z)\left(\omega_1\omega_2(z,z)-|\lambda|^2\right)^{-1}\end{bmatrix}.
\end{multline}
Inverting under the expectation in the left hand side of \eqref{R-dia} and taking the expectation
gives
$$
\begin{bmatrix}
z\phi\left(\left[z^2-(\lambda-x)(\lambda-x)\gwia\right]^{-1}\right)  &
\phi\left((\lambda-x)\left[(\lambda-x)\gwia(\lambda-x)-z^2\right]^{-1}\right)\\
\phi\left((\lambda-x)\gwia\left[(\lambda-x)(\lambda-x)\gwia-z^2\right]^{-1}\right) &
z\phi\left(\left[z^2-(\lambda-x)\gwia(\lambda-x)\right]^{-1}\right)
\end{bmatrix}
$$
Traciality of $\phi$ easily implies that the $(1,1)$ and $(2,2)$ entries of the above matrix are
equal, guaranteeing thus that $\omega_1(z,z)=\omega_2(z,z)=\omega(z)$. Since the above matrix must be equal to the one in \eqref{trente-neuf}, solving a quadratic
equation and recalling that the asymptotics at infinity of $\omega(z)$ is of order $z$ allows us
to write
$$
\omega(z)=\frac{1+\sqrt{1+4z^2|\lambda|^2\phi\left(\left[z^2-(\lambda-x)\gwia(\lambda-x)\right]^{-1}\right)^2}}{2z\phi\left(\left[z^2-(\lambda-x)\gwia(\lambda-x)\right]^{-1}\right)}.
$$
The regularized Cauchy-Stieltjes transform $G_{\mu_x,\epsilon}$ of $x$, namely the $(2,1)$ entry
$\phi\left((\lambda-x)\gwia\left[(\lambda-x)(\lambda-x)\gwia+\epsilon^2\right]^{-1}\right)$,
is then given by
\begin{eqnarray}\label{Gofomega}
G_{\mu_x,\epsilon}(\lambda)&=&\frac{\overline{\lambda}}{|\lambda|^2-\omega(i\epsilon)^2}\\
& =&\frac{\overline{\lambda}}{|\lambda|^2-\frac{\left(1+\sqrt{1-4\epsilon^2|\lambda|^2\phi\left(\left[\epsilon^2+(\lambda-x)\gwia(\lambda-x)\right]^{-1}\right)^2}\right)^2}{2\epsilon^2\phi\left(\left[\epsilon^2+(\lambda-x)\gwia(\lambda-x)\right]^{-1}\right)^2}}.
\nonumber
\end{eqnarray}
Quite trivially, if $\lambda$ does not belong to the spectrum of $x$ and we let $\epsilon$ go to
zero, then  $G_{\mu_x}(\lambda)=\frac{1}{\lambda}$. In particular, this holds for $|\lambda|>
\|x\|$.

Now, a direct computation shows that applying the fixed point equation determining the
subordination functions allows us to write
\begin{eqnarray*}
\omega_1(z,w)&=&\frac{|\lambda^2|}{\omega_2(z,w)}+\frac{1}{\phi\left(\frac{w-\frac{|\lambda^2|}{\omega_1(z,w)}}{\left(z-\frac{|\lambda^2|}{\omega_2(z,w)}\right)\left(w-\frac{|\lambda^2|}{\omega_1(z,w)}\right)-xx\gwia}\right)}\\
\omega_2(z,w)&=&\frac{|\lambda^2|}{\omega_1(z,w)}+\frac{1}{\phi\left(\frac{z-\frac{|\lambda^2|}{\omega_2(z,w)}}{\left(z-\frac{|\lambda^2|}{\omega_2(z,w)}\right)\left(w-\frac{|\lambda^2|}{\omega_1(z,w)}\right)-x\gwia x}\right)}.
\end{eqnarray*}
In particular, for $z=w$ we have seen that $\omega_1(z,z)=\omega_2(z,z)=\omega(z)$, and
\begin{equation}\label{omega-zero}
\omega(z)=\frac{|\lambda^2|}{\omega(z)}+\frac{1}{\phi\left(\frac{z-\frac{|\lambda^2|}{\omega(z)}}{\left(z-\frac{|\lambda^2|}{\omega(z)}\right)^2-x\gwia x}\right)}=
\frac{|\lambda^2|}{\omega(z)}+\frac{1}{\phi\left(\frac{z-\frac{|\lambda^2|}{\omega(z)}}{\left(z-\frac{|\lambda^2|}{\omega(z)}\right)^2-a^2}\right)}
\end{equation}
From this expression and the relation \eqref{trente-neuf} we recognize the (well-known)
fact that the distribution of our $R$-diagonal operator $x$ depends only on its positive part. This functional
equation guarantees at the same time that the dependence of $\omega(z)$ on the {\em argument} of $\lambda$ is constant, thus guaranteeing (via \eqref{Gofomega}) that the distribution of
$x$ has radial symmetry. Now we shall describe the precise dependence of this distribution on
the distance from zero: by simple algebraic manipulations, relation \eqref{omega-zero} becomes
\begin{equation}\label{pre-psi}
\phi\left(\frac{a^2}{\left(z-\frac{|\lambda|^2}{\omega(z)}\right)^2-a^2}\right)=
\frac{z-\omega(z)}{\omega(z)-\frac{|\lambda|^2}{\omega(z)}}
\end{equation}

In terms of the analytic transform $\psi$ we can write \eqref{pre-psi} in the form
\begin{equation}\label{psi}
\psi_{\mu_{a^2}}\left(\left(z-\frac{|\lambda|^2}{\omega(z)}\right)^{-2}\right)=
\frac{z-\omega(z)}{\omega(z)-\frac{|\lambda|^2}{\omega(z)}}.
\end{equation}
We recall that $\omega(i\epsilon)\in i\mathbb R_+$, and in fact $\Im\omega(i\epsilon)>\epsilon
$. We define
$$
A_\epsilon(\lambda)=\Im\frac{1}{i\epsilon-\frac{|\lambda|^2}{\omega(i\epsilon)}}=\frac{1}{
\epsilon+\frac{|\lambda|^2}{\Im\omega(i\epsilon)}},\quad \epsilon>0.
$$
This transforms \eqref{psi} into
\begin{equation}\label{psi2}
\psi_{\mu_{a^2}}\left(-A_\epsilon(\lambda)^2\right)=
\frac{\epsilon A_\epsilon(\lambda)-(\epsilon^2+|\lambda|^2)A_\epsilon(\lambda)^2}{(\epsilon^2+|\lambda|^2)A_\epsilon(\lambda)^2-2\epsilon A_\epsilon(\lambda)+1}.
\end{equation}
Recall from the formula of $\omega(z)$ that
$$
\lim_{\epsilon\to0}\epsilon\Im\omega(i\epsilon)=
\lim_{\epsilon\to0}\frac{1+\sqrt{1-4\epsilon^2|\lambda|^2\phi\left(\left[\epsilon^2+(\lambda-x)\gwia(\lambda-x)\right]^{-1}\right)^2}}{2\phi\left(\left[\epsilon^2+(\lambda-x)\gwia(\lambda-x)\right]^{-1}\right)}.
$$
This quantity (while depending on $\lambda$) is necessarily positive and finite, zero almost everywhere
in the spectrum of $x$. For all these $\lambda$, we take the limit as $\epsilon\to 0$ in \eqref{psi2}
to obtain
$$
\psi_{\mu_{a^2}}\left(-A_0(\lambda)^2\right)=\frac{|\lambda|^2(-A_0(\lambda)^2)}{1
-|\lambda|^2(-A_0(\lambda)^2)}.
$$
We should note that this functional equation is quite trivially solvable on $(-\infty,0)$.
Indeed, with the obvious notations, $\psi_{\mu_{a^2}}(f(r))+1=\frac{1}{1-r^2f(r)}$ is equivalent
to $r^2=\frac{\eta_{\mu_{a^2}}(f(r))}{f(r)}$ (recall that $\eta=\frac{\psi}{1+\psi}$), and $g\colon
v\mapsto\frac{\eta_{\mu_{a^2}}(v)}{v}=\frac{1}{v}-F_{\mu_{a^2}}\left(\frac{1}{v}\right)$ is known
\cite{Belinschi-BercoviciIMRN} to be injective, in fact strictly increasing whenever $a\not\in
\mathbb C\cdot1$, on $(-\infty,0)$ with $0^-\mapsto\phi(a^2)$.
This might be a more convenient way to express $\omega$: we would have
$$
G_{\mu_{x,\epsilon}}(\lambda)=\frac{\overline{\lambda}}{|\lambda|^2\left(1+\frac{|\lambda|^2
A_\epsilon(\lambda)^2}{1+\epsilon A_\epsilon(\lambda)}\right)},
$$
and when $\epsilon A_\epsilon(\lambda)\to0$ as $\epsilon\to0$,
$$
G_{\mu_x}(\lambda)=\frac{\overline{\lambda}}{|\lambda|^2\left(1+|\lambda|^2
A_0(\lambda)^2\right)}
$$
Since
$$
A_0(\lambda)=\sqrt{-g^{\langle-1\rangle}(|\lambda|^2)},
$$
we write
\begin{equation}\label{Gg}
G_{\mu_x}(\lambda)=\frac{1}{\lambda\left(1-|\lambda|^2
g^{\langle-1\rangle}(|\lambda|^2)\right)}.
\end{equation}
The border $|\lambda|^2=\phi(a^2)$ follows from the above remark on the behavior next to
zero of $g$.
However, equally easily in terms of the $S$-transform, we have
$$
S_{\mu_{a^2}}^{\langle-1\rangle}\left(|\lambda|^{-2}\right)=\frac{|\lambda|^2(-A_0(\lambda)^2)}{1
-|\lambda|^2(-A_0(\lambda)^2)}=\frac{\omega(0)^2}{|\lambda|^2-\omega(0)^2}.
$$
Then
\begin{equation}\label{GS}
G_{\mu_x}(\lambda)=\frac{1}{\lambda}\left(1+S_{\mu_{a^2}}^{\langle-1\rangle}\left(|\lambda|^{-2}\right)\right),
\end{equation}
with the exact same restrictions on $\lambda$ as above.

So we have finally
$$ G_{\mu_x}(\lambda)=\begin{cases} \frac{1}{\lambda} & \text{for } |\lambda|\geq \phi(a^2),
\\ \frac{1+S_{\mu_{a^2}}^{\langle -1\rangle}(
|\lambda|^{-2})} {\lambda} & \text{for } |\lambda|^2\leq \phi(a^2),
\end{cases}
$$
which allows us to recover the following theorem of Haagerup and Larsen \cite{HaagerupLarsen}.
\begin{theorem}
The 
Brown measure $\mu_x$ of the $R$-diagonal operator $x=ua$ is the unique rotationally invariant
probability measure such that 
$$\mu_x\big\{\lambda\in
\mathbb{C}:|\lambda|\leq z\big\}= \begin{cases}
0 & \text{for } z\leq \frac{1}{\sqrt{ \phi\left( (xx\gwia)^{-1} \right)  }},  \\
1+S_{xx\gwia}^{\langle -1 \rangle} \left(z^{-2}\right) & \text{for
} \frac{1}{\sqrt{ \phi\left( (xx\gwia)^{-1} \right)  }} \leq z\leq \sqrt{\phi(x x\gwia)}, \\
1 & \text{for } z\geq \sqrt{\phi(x x\gwia)}.
\end{cases} $$

\end{theorem}

Let us recall that Guionnet, Krishnapur, and Zeitouni showed in \cite{GZ} that this Brown measure describes indeed the asymptotic eigenvalue distribution of the corresponding bi-unitarily invariant random matrices.

\section{Brown measure of elliptic--triangular operators}
\label{sec:triangular}

\subsection{Elliptic and triangular--elliptic ensembles}
\label{subsubsec:elliptic} Consider a random matrix
$A_N=(a_{ij})_{1\leq i,j\leq N}$ such that the joint distribution
of random variables $(\Re a_{ij}, \Im a_{ij})_{1\leq i,j\leq N}$
is centered Gaussian with the covariance given by
\begin{align*} \E a_{ij} \overline{a_{kl}}= & \begin{cases} \frac{\alpha}{N}
\delta_{ik} \delta_{jl}  &
\text{if } i<j, \\
\frac{\alpha+\beta}{2 N} \delta_{ik} \delta_{jl}  & \text{if } i=j, \\
\frac{\beta}{N} \delta_{ik} \delta_{jl}  & \text{if } i>j,
\end{cases} \\
\E a_{ij} a_{kl}= & \frac{\gamma}{N} \delta_{il} \delta_{jk},
\end{align*}
where $\alpha,\beta\geq 0$ and $\gamma\in\C$ are such that
$|\gamma|\leq \sqrt{\alpha \beta}$. Informally speaking: there is
a correlation between $a_{ij}$ and $a_{ji}$ and random variables
$(\Re a_{ij},\Im a_{ij})_{1\leq i,j\leq N}$ are as independent as
it is possible to fulfill this requirement. Furthermore, the
entries above the diagonal have the same variance; also the
entries below the diagonal have the same variance (but these two
variances need not to coincide).

We call $A_N$ triangular--elliptic random matrix. It is a natural
generalization of some important random matrix ensembles: for
$\alpha=\beta$ it coincides with the elliptic ensemble and in particular for $\alpha=\beta=1$,
$\gamma=0$ it coincides with the Wigner matrix (i.e.\ $(\Re
a_{ij},\Im a_{ij})$ is a family of iid Gaussian variables) and for
$\alpha=\beta=\gamma=1$ it is a random Hermitian matrix which
coincides with the Gaussian Unitary Ensemble. For $\alpha=1$ and
$\beta=\gamma=0$ the sequence $(A_N)$ converges in
$\star$--moments to a very interesting quasinilpotent operator $T$
(see, e.g., 
\cite{DyHa,Sniady2002,AH})
and for $\alpha=\sqrt{1+t^2}$, $\beta=t$, $\gamma=0$ the sequence
$(A_N)$ converges in $\star$--moments to $T+t Y$, where $Y$ is the
Voiculescu circular element such that $T$ and $Y$ are free. The
Brown measure of the latter operator was computed by Aagaard and
Haagerup \cite{AH}.

One can show
\cite{DyHa} that
the sequence of random matrices $A_N$ converges in
$\star$--moments to a certain generalized circular element $x$
which will be described precisely below in Section
\ref{subsec:triangular}.

\subsection{Elliptic triangular operators}
\label{subsec:triangular}

In this section we will use the notions of operator--valued free
probability; the necessary notions can be found in
\cite{Speicher1998}.

\subsubsection{Preliminaries}

Let $\alpha,\beta\geq 0$, $\gamma\in\C$ such that $|\gamma|\leq
\sqrt{\alpha \beta}$ be fixed. Let $\B=\El^{\infty}(0,1)$, let
$(\B\subset\A,\widetilde{\E}:\A\rightarrow\B)$ be an
operator--valued probability space and let $x\in\A$ be a
generalized circular element $x$ the only nonzero free cumulants
of which are given by
\begin{align}
 \label{eq:kowaria1} k(x, f x)(t)= & \gamma \int_0^1 f(s) ds, \\
 \label{eq:kowaria2} k(x\gwia, f x\gwia)(t)=& \overline{\gamma} \int_0^1 f(s) ds, \\
 k(x, f x\gwia)(t)=& \alpha
\int_t^1 f(s) ds+ \beta \int_0^t f(s) ds, \\
 \label{eq:kowaria4} k(x\gwia, f x)(t)=&
\alpha \int_0^t f(s) ds+ \beta \int_t^1 f(s) ds
\end{align}
for every $f\in\B$. The reader may find the details of this
construction in the case $\alpha=1$, $\beta=\gamma=0$ in
\cite{Sniady2002}.

In order to be able to consider the Brown measure of $x$ we need
to define a tracial state $\phi:\A\rightarrow\C$. We do this by
setting
$$\phi(f)=\int_0^1 f(s) ds $$
for $f\in\B$ and in the general case $\phi(y)=\phi\big(
\widetilde{\E}(y) \big) $. One can show
\cite{DyHa} that
the sequence of random matrices $A_N$ considered in Section
\ref{subsubsec:elliptic} converges in $\star$--moments to $x$ and
therefore $\phi$ is indeed a tracial state.

\subsubsection{Computation of regularized Cauchy transform}

Since we are dealing with an operator--valued case it is useful to
define $\E:\M_2(\A)\rightarrow\M_2(\B)$ as
$$ \E
\begin{bmatrix} a_{11} &
a_{12} \\ a_{21} & a_{22} \end{bmatrix} = \begin{bmatrix} \widetilde{\E}(a_{11}) & \widetilde{\E}(a_{12}) \\
\widetilde{\E}(a_{21}) & \widetilde{\E}(a_{22}) \end{bmatrix} $$
instead of the definition \eqref{eq:conditional}.

The relation between the free cumulants of $x$ and the free
cumulants of $\x$, as defined in \eqref{eq:X}, implies that $X$ is an operator-valued semicircular element whose
$R$--transform is explicitly given by
$$
R_X \begin{bmatrix} a_{11} & a_{12} \\ a_{21} & a_{22}
\end{bmatrix} =
\begin{bmatrix}   k(x, a_{22} x\gwia) & k(x,a_{21} x)  \\
k(x\gwia,a_{12} x\gwia) & k(x\gwia, a_{11} x)
\end{bmatrix}.$$
Thus the general equation 
$G(\Lambda)=(\Lambda-R(G(\Lambda)))^{-1}$
for an operator-valued semicircular element gives in our case
\begin{equation}
\label{eq:warszawa} \G_{\epsilon}(\lambda) =
\begin{bmatrix} g_{11} & g_{12} \\ g_{21} & g_{22} \end{bmatrix}=
 \frac{1}{\ddd}
\begin{bmatrix} i \epsilon- k(x\gwia, g_{11} x) & -\lambda+k(x,g_{21} x) \\
-\bar{\lambda}+k(x\gwia,g_{12} x\gwia) & i \epsilon-k(x,g_{22}
x\gwia)
\end{bmatrix},
\end{equation}
where \begin{multline} \label{eq:det} \ddd=\det \begin{bmatrix}
i\epsilon-k(x,g_{22} x\gwia) & \lambda-k(x,g_{21} x) \\
\bar{\lambda}-k(x\gwia,g_{12} x\gwia) & i\epsilon -k
(x\gwia,g_{11} x)
\end{bmatrix}= \\ \big(i \epsilon-k(x\gwia, g_{11} x)\big)\big(i
\epsilon-k(x,g_{22} x\gwia) \big)-\big(\lambda -k(x,g_{21}
x)\big)\big(\overline{\lambda}-k(x\gwia,g_{12} x\gwia) \big).
\end{multline}
For simplicity, we suppressed the dependence on $\epsilon$ and $\lambda$ of
$g_{ij}=g_{\epsilon,\lambda,ij}\in\B$ and
$\ddd=\ddd_{\epsilon,\lambda}\in\B$.

Observe that for fixed $\epsilon$ and $\lambda$ by
\eqref{eq:kowaria1}--\eqref{eq:kowaria4} the second summand on the
right--hand side of \eqref{eq:det} is a constant function in $\B$
and
\begin{multline} \ddd{}'= \big(i \epsilon-k(x\gwia, g_{11} x)\big)'
\big(i \epsilon- k(x,g_{22} x\gwia) \big)+\\ \shoveright{\big(i
\epsilon-k(x\gwia, g_{11} x)\big)\big(i
\epsilon-k(x,g_{22} x\gwia) \big)'= }\\
-(\alpha-\beta) g_{11} \big(i \epsilon-k(x,g_{22} x\gwia) \big) -
\big(i \epsilon-k(x\gwia, g_{11} x)\big) (\beta-\alpha) g_{22} =0,
\end{multline}
where the last equality follows from the comparison of the matrix
entries in \eqref{eq:warszawa}. It follows that $\ddd \in\B$ is in
fact a constant function. By comparing the entries of
\eqref{eq:warszawa} we see that also $g_{12},g_{21}\in\B$ are
constant. 
Equation \eqref{eq:inwers} implies that $\ddd\in\R$. 

The comparison of upper--left corners of \eqref{eq:warszawa} gives
us a simple integral equation for the function $g_{11}$ which has
a unique solution
\begin{equation} \label{eq:rozwiazaniea} g_{11}(t) = \begin{cases} \frac{i \epsilon
(\alpha-\beta)}{\left( \alpha-\beta e^{\frac{\beta-\alpha}{\ddd}}
\right) \ddd} e^{\frac{\beta-\alpha}{\ddd} t} & \text{for }
\alpha\neq \beta, \\ \frac{i \epsilon}{\ddd+\alpha} & \text{for }
\alpha=\beta. \end{cases} \end{equation}
 In fact, the case $\alpha=\beta$ may be regarded as
a special case of $\alpha\neq \beta$ by taking the limit
$\alpha\to\beta$. Similarly, we find
\begin{equation} \label{eq:rozwiazanied} g_{22}(t) =
\begin{cases} \frac{i \epsilon (\alpha-\beta)}{\left( \alpha-\beta
e^{\frac{\beta-\alpha}{\ddd}} \right) \ddd}
e^{\frac{\beta-\alpha}{\ddd} (1-t)} & \text{for }\alpha\neq \beta,
\\ \frac{i \epsilon}{\ddd+\alpha} & \text{for } \alpha=\beta. \end{cases} \end{equation} We also find
$$g_{12}=\frac{-\gamma \overline{\lambda}-\lambda
\ddd}{\ddd^2-|\gamma|^2}, $$
\begin{equation}
\label{eq:g21} g_{21}=\frac{-\overline{\gamma}
\lambda-\overline{\lambda} \ddd}{\ddd^2-|\gamma|^2}.
\end{equation}

Hence
\begin{equation}\label{eq:rownanie1}
\frac{1}{\ddd}=\det \begin{bmatrix} g_{11} & g_{12} \\ g_{21} &
g_{22} \end{bmatrix} = -\frac{\epsilon^2 (\alpha-\beta)^2
e^{\frac{\beta-\alpha}{\ddd}}}{\left( \alpha-\beta
e^{\frac{\beta-\alpha}{\ddd}} \right)^2 \ddd^2}-\frac{\big| \gamma
\overline{\lambda}+\lambda \ddd\big|^2}{\big(
\ddd^2-|\gamma|^2\big)^2 }
\end{equation} for $\alpha\neq \beta$
and
\begin{equation}\label{eq:rownanie2}
\frac{1}{\ddd}=\det \begin{bmatrix} g_{11} & g_{12} \\ g_{21} &
g_{22} \end{bmatrix}   =
-\frac{\epsilon^2}{(\ddd+\alpha)^2}-\frac{\big| \gamma
\overline{\lambda}+\lambda \ddd\big|^2}{\big(
\ddd^2-|\gamma|^2\big)^2 }\end{equation} if $\alpha=\beta$.

A priori, all above statements hold only for $\epsilon$ in some
neighborhood of infinity, but it is easy to check from the
definition that $g_{\epsilon,\lambda,ij}(s)$ is an analytic
function of $\epsilon$ (in the region $\Re \epsilon\neq 0$) for
fixed values of $0\leq s\leq 1$, $i,j\in\{1,2\}$ and
$\lambda\in\C$. Similarly, $\ddd=(g_{11} g_{22}-g_{12}
g_{21})^{-1}$ is an analytic function of $\epsilon$. It follows
that the above identities hold for all $\epsilon>0$.

\subsubsection{Existence and continuity of Cauchy transform}

Since the function $x\mapsto e^x$ is convex it follows that
$$\frac{\alpha-\beta}{\log \alpha-\log \beta}= \frac{1}{\log
\alpha-\log \beta} \int_{\log \beta}^{\log  \alpha} e^x dx \geq
e^{\frac{1}{\log \alpha-\log \beta} \int^{\log \alpha}_{\log
\beta} x dx}=\sqrt{\alpha \beta}.$$ It is easy to check from
\eqref{eq:inwers} that for fixed $\lambda\in\C$ and
$\epsilon\to\infty$ we have $\ddd\to -\infty$. In the case
$\alpha\neq \beta$ equation \eqref{eq:rozwiazaniea} implies that
$\ddd \neq \frac{\beta-\alpha}{\log \alpha-\log \beta}$ for all
$\epsilon>0$. From the Darboux property it follows that $\ddd <
\frac{\beta-\alpha}{\log \alpha-\log \beta}\leq -|\gamma| $ for
$\epsilon>0$. Similarly, if $\alpha=\beta$ one can show that
$\ddd< -\alpha\leq -|\gamma|$.

Observe that for fixed $\lambda\in\C$ and $\ddd\in\R$ there is at
most one $\epsilon>0$ for which \eqref{eq:rownanie1} holds true,
therefore the continuous function
$\epsilon\mapsto\ddd_{\lambda,\epsilon}$ must be monotone. Since
$\lim_{\epsilon\to\infty} d_{\epsilon,\lambda}=-\infty$ hence
$\epsilon\mapsto\ddd_{\lambda,\epsilon}$ must be decreasing and
the limit $\ddd_{0,\lambda}:=\lim_{\epsilon\to 0}
\ddd_{\epsilon,\lambda}$ exists and is finite. Furthermore, except
for the trivial case $\alpha=\beta=\gamma=0$ we have
$\ddd_{0,\lambda}<0$. Equation \eqref{eq:rownanie1} implies that
$\ddd=\ddd_{0,\lambda}$ is a solution of the equation
\begin{multline}
\label{eq:potwor} \left( \alpha-\beta
e^{\frac{\beta-\alpha}{\ddd}} \right)^2 \left[ \frac{\big(
\ddd^2-|\gamma|^2\big)^2}{\ddd} + \big| \gamma
\overline{\lambda}+\lambda \ddd\big|^2\right] = \\
-\frac{\epsilon^2 (\alpha-\beta)^2 e^{\frac{\beta-\alpha}{\ddd}}
\big( \ddd^2-|\gamma|^2\big)^2}{\ddd^2}
\end{multline}
with $\epsilon=0$.

It is easy to see that for each $\lambda\in\C$ and $\epsilon=0$
there are only finitely many (at most $5$) solutions $\ddd\in\R$ of
the above equation, therefore for every $\lambda_0\in\C$ and
$\delta_0>0$ we can find $0<\delta<\delta_0$ such that the triples
$\lambda_0$, $\epsilon=0$, $\ddd=\ddd_{0,\lambda_0}\pm \delta$ are
not the solutions. It follows that there exists $\epsilon_0>0$
such that for all  $0<\epsilon<\epsilon_0$ and
$|\lambda-\lambda_0|<\epsilon_0$ the triples $\lambda$,
$\epsilon$, $\ddd_{0,\lambda_0}\pm\delta$ are not the solutions.
The function $(\epsilon,\lambda)\mapsto \ddd_{\epsilon,\lambda}$
is continuous for $\epsilon>0$ hence from Darboux property it
follows that for $|\lambda-\lambda_0|<\epsilon_0$ and
$0<\epsilon<\epsilon_0$ we have
$|\ddd_{\epsilon,\lambda}-\ddd_{0,\lambda_0}|<\delta$ hence
$|\ddd_{0,\lambda}-\ddd_{0,\lambda_0}|\leq \delta$. It follows
that the function $\lambda\mapsto\ddd_{0,\lambda}$ is continuous.
Equation \eqref{eq:g21} implies that also the Cauchy transform
$\lambda\mapsto G_x(\lambda)$ is continuous---possibly except for
the case $\alpha=\beta=|\gamma|=-\ddd_{0,\lambda}$.

\subsubsection{Computation of non--regularized Cauchy transform}

Equation \eqref{eq:rozwiazaniea} implies that in the limit as
$\epsilon\to 0$, either $g_{11}$ tends uniformly to zero or
$\alpha-\beta e^{\frac{\beta-\alpha}{\ddd_{\epsilon,\lambda}}}\to
0$ for $\alpha\neq\beta$ or $d\to -\alpha$ for $\alpha=\beta$.

Suppose $g_{11}\to 0$; then also $g_{22}\to 0$ and the comparison
of bottom--left corners of \eqref{eq:warszawa} together with
\eqref{eq:det} gives us
\begin{equation}
\label{eq:besancon} \gamma \big( G_x(\lambda) \big)^2 - \lambda
G_x(\lambda) +1 =0,
\end{equation}
hence
\begin{equation} \label{eq:case3} G_x(\lambda)= \frac{-\lambda\pm
\sqrt{\lambda^2-4 \gamma}}{-2 \gamma}= \frac{2}{\lambda \mp
\sqrt{\lambda^2-4 \gamma}}
\end{equation}
for $\gamma\neq 0$ and $G_x(\lambda)= \lambda^{-1}$ for
$\gamma=0$; note that \eqref{eq:case3} is still valid in the
latter case.

Suppose now that $\alpha-\beta e^{\frac{\beta-\alpha}{\ddd}}\to 0$
if $\alpha\neq\beta$ or $d\to -\alpha$ of $\alpha=\beta$; it
follows that
\begin{equation}
\label{eq:case4} G_x(\lambda)=\frac{-\overline{\lambda} \ddd +
\overline{\gamma} \lambda }{\ddd^2- |\gamma|^2},
\end{equation}
where $d$ is given by
$$ \ddd=\ddd_{0,\lambda} = \begin{cases} \frac{\beta-\alpha}{\log \alpha-\log
\beta} & \text{for } \alpha\neq \beta, \\ -\alpha & \text{for }
\alpha=\beta . \end{cases}
$$

Let us summarize the above discussion: we showed that
$\lambda\mapsto G_x(\lambda)$ is a continuous function given at
each $\lambda$ either by \eqref{eq:case3} or by \eqref{eq:case4}.
These two solutions coincide on the ellipse given by the system of
equations \eqref{eq:besancon}, \eqref{eq:case4}. Therefore on each
of the connected components of the complement of the ellipse the
Cauchy transform is given either by \eqref{eq:case3} or by
\eqref{eq:case4}.

At infinity, the Cauchy transform satisfies
$\lim_{|\lambda|\to\infty} G_x(\lambda)=0$ therefore
\eqref{eq:case3} is the correct choice on the outside of the
ellipse.

On the other hand, for $\lambda=0$ the factor (except, possibly,
for the case $\alpha=\beta=|\gamma|$) $$\frac{\big(
\ddd^2-|\gamma|^2\big)^2}{\ddd} + \big| \gamma
\overline{\lambda}+\lambda \ddd\big|^2$$ is non--zero and
therefore equation \eqref{eq:potwor} implies that \eqref{eq:case4}
is the correct choice of the solution on the inside of the
ellipse.

Having computed the Cauchy transform, we can easily compute the
Brown measure.

\begin{theorem}
The Brown measure of the operator $x$ is the uniform probability
measure on the inner part of the ellipse given by the system of
equations \eqref{eq:besancon}, \eqref{eq:case3}.
\end{theorem}

A careful reader might object that for the case $\alpha=\beta=e^{2
i\tau} \gamma$ with $\tau\in\R$ our proof has a gap since we
cannot guarantee that the Cauchy transform $G_x(\lambda)$ is
continuous. However, in this case the operator $e^{-i \tau} x$ is
self-adjoint and coincides with Voiculescu's semicircular operator, and
the calculation of its spectral measure is trivial. On the other hand
one can easily check that in the case $\alpha=\beta=|\gamma|$ our
ellipse degenerates to an interval and the uniform measure on such
a degenerated ellipse coincides with the semicircular measure.
Therefore our result is true also in this case.
That the Brown measure of the operator $x$ is, at least for $\beta\not=0$, indeed also the asymptotic
eigenvalue distribution of the elliptic-triangular random matrices $A_N$ follows from the result of \'Sniady \cite{Sn} that a small Gaussian deformation
of a random matrix ensemble does not change the Brown measure in the limit, but makes the convergence of Brown measure continous; in our case,
a small Gaussian deformation does not change the nature of the considered
ensemble.

\section{Final remarks: Discontinuity of Brown spectral
measure} \label{sec:discontinuity}

Following the program from
Section \ref{subsec:convergence} we would like to find the
connection between the eigenvalue density of the random matrices
$A_N$ and the Brown measure $\mu_x$ of their limit; in particular
one would hope that the Brown spectral measure is continuous with
respect to the topology of the convergence of $\star$--moments and
hence the empirical eigenvalue distributions $\mu_{A_N}$ converge
to $\mu_x$. Unfortunately, in general this is not true; a very
simple counterexample is presented in
\cite{Haagerup2001,Sn}. The reason for this phenomenon is
that the definition of the Fuglede--Kadison determinant uses the
logarithm, a function unbounded from below on any neighborhood of
zero.

Let us consider some sequence $(A_N)$ of random matrices which
converges in $\star$--moments to some $x$. Even though there exist
such sequences with a property that the eigenvalue densities
$\mu_{A_N}$ do not converge to $\mu_x$, there is a growing
evidence that such examples are very rare. In particular, it was
shown by Haagerup \cite{Haagerup2001} and later by \'Sniady
\cite{Sn} that every such sequence can be perturbed by a
certain small random correction in such a way that the new
sequence $(A_N')$ still converges to $x$ and furthermore the Brown
measures converge: $\mu_{A_N'}\to\mu_x$. In general, adding a small
perturbation changes the nature of the considered random matrix model,
so these results do not apply directly to the original ensemble.
(An exception from this is the case of elliptic-triangular random matrices,
considered in Section \ref{sec:triangular}.)

There has been a lot of research in this context in the last ten years;
in particular, controlling the discontinuity of the Brown measure and
thus showing that the asymptotic eigenvalue distribution of the random matrices is indeed given by the Brown measure of the limit operator was achieved in considerable generality in the circular law for Wigner matrices and then, more generally, also for $R$-diagonal operators and bi-unitarily random matrices in \cite{GZ}. At the moment it is not clear whether the ideas
from those investigations apply also to situations like polynomials in Gaussian
(or even Wigner matrices). However, we believe the following conjecture to be true.
We will address this question in forthcoming investigations.

\begin{conjecture}
Let $p$ be a (not necessarily self-adjoint) polynomial in $m$ non-commuting
variables. Consider $m$ independent self-adjoint Gaussian (or, more general, Wigner)
random matrices $X_N^{(1)},\dots, X_N^{(m)}$. One knows that they converge
in $\star$-moments to a free semicircular family $s_1,\dots,s_m$. Consider now the polynomial $p$ evaluated in the random matrices and in the semicircular family, respectively; i.e.,
$$A_N:=p(X_N^{(1)},\dots,X_N^{(m)}) \qquad\text{and}\qquad
x=p(s_1,\dots,s_m).$$
The convergence in $\star$-moments of $X_N^{(1)},\dots, X_N^{(m)}$ to $s_1,\dots,s_m$ implies then that also the (in general, non-normal) $A_N$ converge in $\star$-moments to the operator $x$.
We conjecture that also the eigenvalue distributions $\mu_{A_N}$ of the random matrices $A_N$ converge to the Brown measure $\mu_x$ of the
limit operator $x$.  
\end{conjecture}

It would be interesting to see whether methods from \cite{ORSV}, where the special case of a product of independent elliptic random matrices is treated,
can be extended to more general polynomials. 

\section{Acknowledgments}

For part of the duration of the work on this project, S.T.B. has been supported by a Discovery Grant from NSERC and by an 
Alexander von Humboldt Research Fellowship. 

Research of P.\'S. supported by \emph{Narodowe Centrum Nauki}, grant number  2014/15/B/ST1/00064.

R.S. has been supported by the ERC Advanced Grant ``Non-commutative distributions in free probability'' (grant no. 339760).

We thank Tobias Mai for several helpful discussions in the context of this work.
We also thank the referee for many suggestions which improved the readability.

\bibliographystyle{alpha}

\end{document}